\documentclass[reqno]{amsart}
\usepackage[utf8]{inputenc}
\usepackage{enumitem}
\usepackage{mathdots}
\usepackage{amsmath,amssymb,amsbsy,amsfonts,amsthm,latexsym,amsopn,amstext,mathtools,leftidx,
            amsxtra,euscript,amscd,verbatim,epsf,colordvi,graphics}
\usepackage{ytableau}
\usepackage{color}
\usepackage{pgf,tikz}
\usetikzlibrary{decorations.pathreplacing}
\usetikzlibrary{arrows}
\usepackage{mathrsfs}
\usepackage{graphicx}

\usepackage{hyperref}
\usepackage{cleveref}
\usepackage[margin=11pt,font=footnotesize,labelfont=bf]{caption}
\newtheoremstyle{tplain}{3pt}{3pt}{\rmfamily}{}{\bfseries}{.}{0.5em}{}
\theoremstyle{tplain}

\usepackage{pgf,tikz}
\usepackage{tcolorbox}
\usepackage[enableskew]{youngtab}
\definecolor{darkgreen}{cmyk}{1,0,1,0}
\newtheorem{thm}{Theorem}
\newtheorem{lem}{Lemma}
\newtheorem{ex}{Example}
\newtheorem{cor}{Corollary}
\newtheorem{prop}{Proposition}
\newtheorem{obs}{Remark}
\newtheorem{defi}{Definition}

\def \rm {\mathrm}
\def \mbf {\mathbf}

\def\Z { \mathbb{Z}}

\def\Z { \mathbb{Z}}

\def\N { \mathbb{N}}
\def\LRAII {\mathsf{LR}}
\def \r {\mathbf{r}}
\def \t {\mathbf{t}}
\def \redu {\mathrm{red}}
\def\c { \mathrm{c}}
\def \k {\mathfrak{k}}

\def\shape { \mathrm{shape}}
\def\rem { \mathrm{rem}}
\def\suc { \mathrm{suc}}
\def\a { \mathbf{a}}
\def\b { \mathbf{b}}
\def\s { \mathbf{s}}
\def\B { \mathbf{B}}

\newcommand{\oo}{\color{blue}}

\makeatletter
\newcommand*\bigcdot{\mathpalette\bigcdot@{.5}}
\newcommand*\bigcdot@[2]{\mathbin{\vcenter{\hbox{\scalebox{#2}{$\m@th#1\bullet$}}}}}
\makeatother

\newcommand*\circled[1]{\tikz[baseline=(char.base)]{
            \node[shape=circle,draw,inner sep=1pt] (char) {#1};}}

\newcommand{\YT}[3]{
\vcenter{\hbox{
\begin{tikzpicture}[x={(0in,-#1)},y={(#1,0in)}] 
\foreach \rowi [count=\i] in {#3} {
 \foreach \e [count=\j] in \rowi {
  \draw (\i,\j) rectangle +(-1,-1);
  \draw (\i-0.5,\j-0.5) node {$#2\e$};
 }
}
\end{tikzpicture}
}}
}

\title[Reverse Schensted bumping, inverse reduction map and concatenation] {Reverse Schensted bumping, inverse reduction map and
concatenation in the inverse \\quantum Littlewood-Richardson bijection}

\author{Olga Azenhas}
\address{ University of Coimbra, CMUC, Department of Mathematics, Portugal}
\email{oazenhas@mat.uc.pt}



\begin{document}

\begin{abstract} Provided the explicit formula  for the inverse reduction map,  arXiv:2606.24840, we exhibit for  symplectic bumped columns disjoint of the action of the parity involution, and for a family of bumped  columns from symplectic intervals closed for the parity  involution, the feasibility of the concatenation operation    to return a semi-standard  tableau in  the inverse  quantum Littlewood-Richardson
bijection. This analysis, includes symplectic bumped columns of the form $x\B y$ such that $\B$ is an interval closed for the parity involution with $x$ and $y$ single entries. 
\end{abstract}
\maketitle

\tableofcontents

\section{Introduction}
{The quantum Littlewood-Richardson bijection has been introduced by Watanabe in \cite{watanabe}.
 The purpose of this note is to make precise  the feasibility of the concatenation operation, given the explicit formula for the inverse
reduction map in \cite{azreduction}, in the inverse  quantum Littlewood-Richardson
bijection to return a semi-standard  tableau.
More precisely, we analyse when the action of the inverse reduction map, $\redu^{-1}$,
on a symplectic reverse bumped column $\b$ from a semi-standard tableau $S\in SST_{2n}(\mu)$ 
 returns a column $\redu^{-1}(\b)$ whose concatenation with the returned reverse bumped semi-standard tableau  $(  \longleftarrow S)$, via the
 reverse Schensted column bumping on $S$,  is a semi-standard tableau. This analysis is carried out in Section \ref{sec:essential} and it has been started in \cite[Section 4.3, Section 4.6, Theorem 9]{azreco} that we here extend to symplectic bumped columns from symplectic intervals closed for the parity involution.

Recall that the parity involution $\s$, Definition \ref{s}, swaps an even number with the previous odd number and an odd number with the next even number. In \cite{azreduction} 
we have shown that
any symplectic column ${\bf a}\in SpT_{2n}(\varpi_t)$ can be decomposed, for some $k\ge 0$, in the form
\begin{align}\label{introdu:symplecticolumn}
{\bf a}=
X^{(0)}{\bf A_1} X^{(1)}{\bf A_2}\cdots {\bf A_{k-1}}X^{(k-1)} {\bf A_k} X^{(k)},
\end{align}

where
\begin{itemize}
\item  for $i=0,1,\dots,k$, each $X^{(i)}$ is a column word of length $p_i\ge 0$, such that

\begin{align}\label{introdu:xx}\s(X^{(0)}X^{(1)}\cdots X^{(k)})\cap {\bf a}=\emptyset, \mbox{ and}
\end{align}

\item for $i=1,\dots, k$, each factor ${\bf A_i}$ is
a maximal nonempty interval  of even length $t_i\ge 2$     such that
\begin{align}\s({\bf A_i})={\bf A_i},
\end{align}
with $$t=\sum_{i=1}^kt_i +\sum_{i=0}^k p_i  \in [n].$$
\end{itemize}

For $i=1,\dots,k$, such intervals are symplectic pieces of the form
$${\bf A_i}=(a_i, a_i+1,\dots,a_i+t_i-1), \mbox{ $a_i\notin 2\Z$}$$
satisfying
\begin{align}
&a_i\ge 2t^{(i-1)}+t_i+1, ~\mbox{ with $\displaystyle t^{(i-1)}:=\sum_{j=1}^{i-1}t_j+\sum_{j=0}^{i-1}p_j$, $t^{(0)}:=p_0$, and $t^{(k)}:=t$.}
   \end{align}

Our current contribution for the aforesaid problem comprises  bumped  columns from symplectic intervals closed for the parity  involution, namely theorems \ref{thm:reducbumpinterval1}, \ref{thm:reducbumpinterval1m+}  and \ref{thm:bumpedxintervaly}, in Subsection  \ref{bumpintervalsymplecticinterval},   and  bumped symplectic intervals from semi-standard tableaux in Theorem \ref{thm:reducbumpintervalgeneral},  Subsection \ref{bumpedinterval}, and for symplectic bumped columns disjoint of the action of the parity involution in Theorem \ref{thm:essentialss}, Subsection  \ref{bumpedcolumndisjoint}. The analysis of bumped columns from symplectic intervals closed for the parity involution is essential due to the canonical form of a symplectic column as shown above in \eqref{introdu:symplecticolumn}.

\section*{Acknowledgements}
The author acknowledges financial support by the Centre for Mathematics of the University of Coimbra (CMUC,
https://doi.org/10.54499/UID/00324/2025) under the Portuguese Foundation for Science and Technology (FCT), Grants UID/00324/2025 and
UID/PRR/00324/2025.

\section{Preliminaries} \label{sec:pre}
\subsection{Partitions and semistandard-tableaux} A \emph{partition} $\gamma$ is a weakly decreasing sequence of nonnegative integers $\gamma_1\ge \gamma_2 \ge \cdots$
such that $\gamma_k =0$ for some $k \ge 1$. The   maximal $i$ such that $\gamma_i> 0$ is called the \emph{number}
of \emph{parts} or \emph{length} of $\gamma$, denoted $\ell(\gamma)$. For each $m\ge 0$, the set of partitions of length at most $m$ is
denoted by \emph{$Par_{\le m}$}. We assume the inclusion $Par_{\le m}\subseteq Par_{\le k}$ whenever $k\ge m$. Thus we often write the
partition $\gamma$ as a vector
 $\gamma = (\gamma_1, \gamma_2, \dots,\gamma_k)$ for $k \ge  \ell(\gamma)$. The \emph{empty} \emph{partition} is the empty sequence $()$ and
 is regarded as the unique partition of length zero.


 A partition $\gamma$ is identified
with its \emph{Young diagram} $D(\gamma)$ which is a left and top justified collection of boxes (or cells)
with $\gamma_k$ many boxes in the $k$th row for all $k \in \mathbb{Z}_> 0$. In particular, the empty Young diagram and  the  partition $()$
are identified. The number of cells of $D(\gamma)$ is  the sum of the parts of $\gamma$ and is denoted by $|\gamma|$.
The boxes or cells of the Young diagram of $\gamma$ are identified by its  coordinates $(i,j)$ in the matrix style, that is, $1\le i\le
\ell(\gamma)$ and $1\le j\le \gamma_i$.

Let $\gamma$, $\mu$  partitions   with $\mu \subseteq \gamma$, that is,  $\mu_i\le  \gamma_i$ for all $ i\in\mathbb{Z}_>0 $, or the Young
diagram of $\mu$ is a subset of
the Young diagram of $\gamma$. The skew-diagram (or skew Young diagram) $\gamma/\mu$ is defined to be $D(\gamma)\setminus D(\mu)$.
For $\lambda\subseteq\gamma$ partitions, we write $\lambda\subseteq_{vert} \gamma$ to mean that  $\gamma/\lambda$  is a \emph{vertical}
\emph{strip}, that is, the skew-diagram $\gamma/\lambda$ has at most one box in each row. The number of cells of $\gamma/\lambda$ is
$|\gamma/\lambda|:=|\gamma|-|\lambda|.$

 A \emph{
tableau} $T$ of (skew) shape $\gamma/\mu$ is a map (or a filling of $D(\gamma)$)
 $$T:D(\gamma)\rightarrow \mathbb{Z}_\ge0, \; (i,j)\mapsto T(i,j),$$
assigning a positive integer to each box of $D(\gamma)\setminus D(\mu)$ and   $0$  to the boxes of $D(\mu)$. Denote by $\rm{Tab}(\gamma/\mu)$
the set of tableaux of shape $\gamma/\mu$.  We say that the tableau $T$ is \emph{semi-standard} if an addition the assignment is
such that it is weakly increasing as we go from left to right along a row and strictly
increasing as we go from top to bottom along a column excluding the boxes in $\mu$,
\begin{align}&T(i,j)\le T(i,j+1),\; T(i,j)<T(i+1,j), \mbox{ for all $(i,j)\in D(\gamma)\setminus D(\mu)$,}\nonumber\\
&\mbox{ and }  T(i,j)=0, \mbox{ for $(i,j)\in  D(\mu)$,}
\nonumber
\end{align}
where we set $T(a,b):=\infty$ if $(a, b) \notin D(\gamma)$.
 Usually $T(i,j)$ is just referred as
the entry in the box $(i,j)$ and we omit the zeroes  in the boxes of $\mu$. A positive integer
$m\ge \ell(\gamma)$ will be fixed and $[0, m]:=\{0, 1, \dots, m\}$ will be used as a co-domain for the map $T$. We call
$[m]:=\{1, \dots, m\}$ the alphabet of the semi-standard tableau $T$.
 In this case, we will denote the \emph{set} of \emph{semi-standard} \emph{tableaux} of \emph{shape} $\gamma/\mu$ by
$SST_m(\gamma/\mu)$. When $\mu=()$, we just write $SST_m(\gamma)$. The \emph{weight} or \emph{content} of $T$ is the nonnegative vector
$\mathrm{wt}(T)=(T[1],\dots,  T[m])$, where $T[i]:=\#\{(a,b)\in D(\gamma):T(a,b)=i\}$ for $i\in[m]$, that is, $T[i]$ is  the number of
occurrences of $i$ in the tableau $T$.
\subsection{Schensted column insertion and  reverse }\label{subsec:insert}  We follow the notation in \cite{azreco}. The Schensted \emph{column} \emph{insertion} or
\emph{column} \emph{bumping} \cite{fulton,stanley},  takes a positive integer $x$ and a tableau $T\in SST(\gamma)$ (semistandard tableau of
shape $\gamma$) and
puts $x$ in a new box at the bottom of the first column (counting from left to right)  if  it is strictly
larger than all the entries of the column. If not, it \emph{bumps} the  \emph{smallest}
\emph{entry}
in the \emph{column} that is \emph{larger} \emph{than} or \emph{equal} to $x$. The bumped entry moves to the
next column, going to the end if possible, and bumping an element to the next column
otherwise. The process terminates when  the bumped entry  goes to the bottom of the next
column, or until it becomes the only entry of a new column. The returned tableau is
denoted $x\bigcdot T:=x \rightarrow T$, that is, $x\bigcdot T$ is the result of column inserting $x$ into $T$, $x \rightarrow T$. We call to $
x \bigcdot T$ the multiplication of $x$ and $T$ in the plactic monoid.

The \emph{reverse Schensted column insertion } takes a tableau $U\in SST(\gamma)$ and an entry $y$ at the bottom of a column of $U$, whose
cell has  \emph{row} \emph{coordinate}, say $r_y$,  and bumps it while deletes its box. If the column is the first, the new tableau is
obtained by deleting the bottom box of the first column of $U$. If not, the bumped $y$ moves  to the next column on its left and \emph{bumps}
the  \emph{largest}
\emph{entry} in the \emph{column} that is \emph{smaller} than or \emph{equal} to $y$. In both cases, if the bumped entry belongs to the first
column the process ends and returns a pair $(b, ~\underset{ r_y}\leftarrow U)$ consisting of the bumped entry $b\le y$ from the first column
of $U$ together with a new  tableau denoted $\underset{r_y}\leftarrow U$, with one box less than $U$. If not, the process continues until an
entry is bumped from the first column of $U$.
\begin{obs}\label{noexplicit} When we want to emphasize that in the reverse Schensted column insertion we are considering  the pair $(b,
~\underset{ r_y}\leftarrow U)$ without making explicit the bumped entry $b$ from the first column of $U$, we just write
\begin{align*}(\quad , ~\underset{ r_y}\leftarrow U).
\end{align*}
\end{obs}

We now recall the bumping routes defined by the column Schensted insertion and by the reverse column Schensted insertion.
\begin{obs}\label{routes}\cite[Appendix A.2, Exercise 3]{fulton}  Let us consider two successive column-insertions of $ x < x'$ into a tableau
$T$. First column-inserting $x$ in $T$ and then column-inserting $x'$ in
the resulting tableau $x\rightarrow T$, yielding  to two \emph{bumping} \emph{routes} $R$ and $R'$, and
two new boxes $B$ and $B'$. Then $R'$ lies strictly below $R$,
and $B'$ is Southwest of $ B$.

If $x\ge x'$, $R'$ lies weakly above $R$ and $B'$ Northeast of $B$.

\end{obs}
Since column insertion is reversible, one has the following \emph{reverse bumping route} property for \emph{reverse} \emph{column}
\emph{insertion}.

\begin{obs}\label{revroutes}
\begin{enumerate}
\item[(a)]
Let $ z$ and   $z'$  be two bottom column entries of a tableau $U$, with $z'$  \emph{strictly} \emph{Southwest} of $z$ (not in the same
row), or let $z<z'$ be the two bottom entries of a column of  $U$. Let $r_z<{r_{z'}}$ be the corresponding cell row coordinates.

 First consider the reverse column-inserting of $z'$  and then reverse column-inserting of $z$ in
the resulting tableau $(\underset{r_{z'}}\leftarrow U)$. This yields   two reverse bumping routes $Z'$ respectively $Z$ such that $Z'$
\emph{lies} \emph{strictly} \emph{below} $Z$.  In particular, if the bumping routes terminating rows terminating rows are respectively
$i'>i$, the
two new entries $a'\le z'$ respectively $a\le z$,   in the first column of $U$ are such that $(\underset{r_{z'}}\leftarrow U)(i',1)=a'
> (\underset{r_{z}}\leftarrow(\underset{r_{z'}}\leftarrow U))(i,1)=a
$ and $i'>i$. The bumped the pair $U(i,1)=b<U(i',1)=b'$ from the first column of $U$ are such that $ (b,b')\le (a,a')\le (z,z')$, that is,
$b\le a\le z$, $b'\le a'\le z'$.
Generalizing the notation in Remark \ref{noexplicit}, we write

\begin{align}(\quad , \underset{ (r_z,r_{z'})}\leftarrow U):=((b,b'),\underset{r_{z}}\leftarrow(\underset{r_{z'}}\leftarrow U))
\end{align}
where  the  entry in the cell $(u,1)$ of the first column of $\underset{ (r_z,r_{z'})}\leftarrow U$ is such that $$ (\underset{
(r_z,r_{z'})}\leftarrow U)(u,1)=\begin{cases}(\underset{r_{z}}\leftarrow(\underset{r_{z'}}\leftarrow U))(i,1)=a& \mbox{ if $u=i$},\\
(\underset{r_{z'}}\leftarrow U)(i',1)=a'& \mbox{ if $u=i'>i$},\\
U(u,1)&  \mbox{ if $u\neq i, i'$}.
\end{cases}$$

 \item[(b)] In case $y'$  \emph{strictly} \emph{Southwest} of $y$ (not in the same row) and we first consider reverse column-inserting of
     $y$  followed with the  reverse column-inserting $y'$ in
the resulting tableau $(\underset{r_{y}}\leftarrow T)$ then the corresponding reverse bumping routes $Z$ and $Z'$ are such that $Z'$
\emph{lies} \emph{weakly} \emph{below} $Z$.
The two new entries $b$ respectively $b'$ in the first column of $T$, $(\underset{r_{y'}}\leftarrow(\underset{r_{y}}\leftarrow
T))(i',1)=b'>(\underset{r_{y}}\leftarrow T)(i,1)=b$ with $i'>i$, or $i=i'$ and $(\underset{r_{y}}\leftarrow T)(i,1)=b\le
(\underset{r_{y'}}\leftarrow(\underset{r_{y}}\leftarrow T))(i,1)=b'$.

 \item [(c)] If two or more columns have the \emph{bottom} \emph{entries} in a \emph{same} \emph{row}, let $y'$ be the \emph{right}
     \emph{most} \emph{bottom} \emph{entry} in that \emph{row}.  Consider the reverse column-inserting of $y'$ in $T$ and then again the
     reverse column-inserting of $y'$ in
the resulting tableau $\underset{r_{y'}}\leftarrow T$. This yields   two bumping routes $Z_1$ respectively $Z_2$ such that $Z_2$ \emph{lies}
\emph{weakly} \emph{below} $Z_1$ and two new entries $b_1$ respectively $b_2$ in the first column of $T$, such that
$(\underset{r_{y'}}\leftarrow(\underset{r_{y'}}\leftarrow T))(i_2,1)=b_2> (\underset{r_{y'}}\leftarrow T))(i_1,1)=b_1$ with $i_2>i_1$, or
$i_1=i_2$ and $(\underset{r_{y'}}\leftarrow T)(i_1,1)=b_1\le (\underset{r_{y'}}\leftarrow(\underset{r_{y'}}\leftarrow T))(i_1,1)=b_2$.
\end{enumerate}
\end{obs}
\medskip

\subsection{Symplectic tableaux}
From now on we fix $n\in\mathbb{N}$.
A semistandard tableau $G \in SST_{2n}(\lambda)$ is said to
be symplectic \cite{king76} if
$$G(k, 1) \ge 2k-1, \mbox{  for all $  k \in  [ \ell(\lambda)]$}.$$
Let $SpT_{2n}(\lambda)$ denote the set of all symplectic tableaux of shape $\gamma$ on the alphabet $[2n]$.


\begin{prop}\cite{azreco,watanabe} \label{prop:symplectic1}  
Let $G \in SST_{2n}(\lambda)$.
\begin{enumerate}
\item\cite{watanabe} If $SpT_{2n}(\lambda)\neq\emptyset$ then $\ell(\lambda)\le n$.
\item If $G=(a_1,\dots,a_l)\in SpT_{2n}(\varpi_l)$ then $l\le n$ and
$l\le \lfloor\frac{a_l+1}{2}\rfloor$.




\item  If $G'$ is obtained from $G\in  SpT_{2n}(\varpi_l)$ by suppressing
$0\le t_0\le l$ entries then $G'\in SpT_{2n}(\varpi_{l-t_0})$ and the suppressed part $G''\in SpT_{2n}(\varpi_{t_0})$.

\item  Let $\mu\in Par_{\le n}$. If $G'$ is obtained from $G\in  SpT_{2n}(\mu)$ by reverse Schensted column insertion, then $G'$ is
    still symplectic.
 \end{enumerate}
\end{prop}

\subsection{The  reduction map in the  quantum Littlewood-Richardson map}\label{sec:containmentmain} For the reader convenience this section recalls several properties of the removal and  reduction maps
in \cite{watanabe}, and on the inverse reduction map in \cite{azreco,azreduction}.
We fix $l \in [0, 2n]$ and $\mathbf{a} = (a_1, \dots , a_l)$ a column  in $SST_{2n}(\varpi_l)$.
We often regard $\mathbf{a}$ as a set.

The \emph{removal} \emph{subword }of $\mathbf{a}$ \cite{watanabe} is defined to be the subword $\mathrm{rem}(\mathbf{a})$ of  $\mathbf{a}$
obtained by the following recursive formula:

\begin{align}\label{removal}
\mathrm{rem}(\mathbf{a}) :=
\begin{cases}
\emptyset,& \mbox{ if } l\le 1,\\
\mathrm{rem}(a_l,\dots, a_{l-2})  (a_{l-1},a_l),& \mbox{ if } l\ge 2, a_l\in   2\mathbb{Z}, a_{l-1} = a_l-1,\\
 &\mbox{ and }
a_l < 2l -|\rem(a_1,\dots,a_{l-2})|-1,\\
\rem(a_1,\dots,a_{l-1})& \mbox{ otherwise}.
\end{cases}
\end{align}

\begin{defi}\cite{watanabe} \label{def:rem} For the column
$\mathbf{a}$, the new column $\mathrm{red}(\mathbf{a})$, reduction of $\mathbf{a}$, is defined to be the one obtained from $\mathbf{a}$ by
removing the entries in
the set $\mathrm{rem}(\mathbf{a})$,
$$\redu( \mathbf{a})=\mathbf{a}\setminus \rem(\mathbf{a})\in SST_{2n}(\varpi_{l-|\rem(\mathbf{a})|}).$$
\end{defi}

\begin{defi}\label{s}For each $x\in\mathbb{Z}$, set \cite{watanabe} \begin{align} \s(x)=\begin{cases}
x+1, &  \mbox{ if } x\notin 2\mathbb{Z}, \\
x-1, &  \mbox{ if }  x\in 2\mathbb{Z}.
\end{cases}
\end{align}
\end{defi}
Indeed $\s^2(x)=x$, for each $x\in \Z$. We call it the parity swapping involution which has the following properties
\begin{align}
min(a_i,\s(a_i))\notin 2\Z,~~min(a_i,\s(a_i))-1\in 2\Z, \\
 max(a_i,\s(a_i))\in 2\Z, ~~ max(a_i,\s(a_i))+1\notin 2\Z.
\end{align}

\begin{prop}\cite[Proposition 4.2.2, 4.2.7, Corollary 4.2.8]{watanabe}\label{prop:rem}
\begin{enumerate}

\item[(1)] If $i\in [1,l]$ and $a_j\notin \rem(\mbf{a})$, for $j\in[i,l]$ then $\rem(\mbf{a})=\rem(a_1,\dots,a_{i-1})$.

\item[(2)] for each $i\in [1,l]$, $a_i\in \rem(\mbf{a})$ if and only if one of the following holds
\begin{enumerate}
\item[(a)] $a_i$ odd, $i<l$, $a_{i+1}=a_i+1$ and $a_i<2i-|\rem(a_1,\dots,a_{i-1})|$

\item[(b)] $a_i$ even, $i>1$, $a_{i}=a_{i-1}+1$ and $a_i<2i-|\rem(a_1,\dots,a_{i-2})|-1$
\end{enumerate}

\item [(3)]  $a_i \in \rem(\mbf{a})$ if and only if $s(a_i) \in \rem(\mbf{a})$.
Consequently, $|\rem(\mbf{a})|\in 2\mathbb{Z}$.


\item[(4)] For $i\in[1,l]$,
\begin{enumerate}
\item[(a)] $a_i$ odd, then $|\rem(a_1,\dots,a_{i-1})|\ge 2i-a_i-1$

\item[(b)] $a_i$ even, then $|\rem(a_1,\dots,a_{i})|\ge 2i-a_i$
\end{enumerate}
\end{enumerate}

\end{prop}

Some properties of the reduction and successor maps   follow. More precisely, the reduction map on columns is the injective assignment as
shown in the next proposition.

\begin{prop} \label{prop:redu}\cite{watanabe,nsw} Let $l\in [0,2n]$ and  $\mathbf{a}=(a_1,\dots,a_l) \in SST_{2n}(\varpi_l)$.
Then
\begin{enumerate}

\item[(1)] $\redu(\mathbf{a})$ is symplectic.

    \item[(2)] $\redu(\mathbf{a})=\mathbf{a}$ if and only if $\mathbf{a}$ is symplectic.

\item[(3)] \cite[Proposition 4.3.6]{watanabe} The reduction map $\redu$ on $SST_{2n}(\varpi_l)$ is an injective assignment

\begin{align}\redu:SST_{2n}(\varpi_l)&\rightarrow\bigsqcup_{\begin{smallmatrix} 0\le t\le min\{l,2n-l\}\\
l-t\in 2\mathbb{Z}
\end{smallmatrix}}SpT_{2n}(\varpi_t)\nonumber\\
\mathbf{a}&\mapsto \redu( \mathbf{a})=\mathbf{a}\setminus \rem(\mathbf{a}).
\label{redumap}
\end{align}
\end{enumerate}
\end{prop}

\subsection{The inverse reduction map in the inverse quantum Littlewood-Richardson bijection}\label{sec:inversereduction}

We now recall results on  the inverse reduction map \eqref{redumap} from
\cite{azreduction} where a complete explicit formula for the inverse reduction map is given. For a fixed $n\in\mathbb{N}$, let $l\in [0,2n]$
and $t\in[0,n]$ such that  $0\le t\le min\{l,2n-l\}$ and $l-t\in 2\Z$. Then we define the expanding map

\begin{align}\label{exp}
\rm{exp}_t:SpT_{2n}(\varpi_t)&\rightarrow SST_{2n}(\varpi_l)\quad\\
\mathbf{a}&\mapsto\rm{exp}(\mathbf{a}):=T,\nonumber
\end{align}
such that $\redu(\rm{exp}(\bf a))=a$. That is, $\rm{red}_{t}^{-1}:=\redu^{-1}_{|SpT_{2n}(\varpi_t)}=\rm{exp}_t$.

Next lemma characterizes the symplectic columns when
the symplectic column is such that  consecutive integers may occur only as an even number followed with an odd number. 

\begin{lem}\label{lem:swapping}\cite{azreduction}Let $\mathbf{a} = (a_1, \dots , a_t)\in SST_{2n}(\varpi_t)$. The following are pairwise
equivalent

\begin{enumerate}
\item for every $1\le i\le t$, $\s(a_i)\notin \bf{a}$
\item ${\bf{a}}\cap \s(\bf{a})=\emptyset$
 \item  $[x, y]$ is an interval factor $\subseteq \bf a $ $\Rightarrow$ $x=y$  $ \vee$ ($x\in 2\Z, y=x+1\notin 2\Z$).
 \item for every $1\le i<t$,
$a_i\notin 2\Z\Rightarrow a_{i+1}-a_i\ge 2$ ( $a_{i+1}\ge \s(a)+1$)
\item for every $1< i\le t$,
$a_i\in 2\Z\Rightarrow a_{i}-a_{i-1}\ge 2$ ( $\s(a)\ge a_{i-1}+ 1$)
\item for every $1\le i<t$,
\begin{align*}
&(  a_i\notin 2\mathbb{Z}, ~~ a_{i+1}\in 2\mathbb{Z}\Rightarrow a_{i+1}-a_i\ge 3)\\
&\qquad\qquad\wedge\\
&( a_i,~ ~a_{i+1}\notin 2\mathbb{Z} \mbox{ or } a_i, ~~a_{i+1}\in 2\mathbb{Z}\Rightarrow a_{i+1}-a_i\ge 2)\\
 &\qquad\qquad\wedge\\
 &( a_i\in 2\mathbb{Z},~ ~a_{i+1}\notin 2\mathbb{Z}\Rightarrow a_{i+1}-a_i\ge 1).
 \end{align*}
\end{enumerate}

\end{lem}

\begin{thm}\cite{ azreduction}\label{thm:nofactors} Let $l\in [0,2n]$ and $t\in[0,n]$ such that  $0\le t\le min\{l,2n-l\}$ and $l-t\in 2\Z$.
Let $\mathbf{a} =
(a_1, \dots , a_t)\in SpT_{2n}(\varpi_t)$ such that $\s(\a)\cap \a=\a$, that is, for every $1\le i\le t$, $\s(a_i)\notin \bf{a}$. Then,

\begin{align}\label{redu:nofactors}
\rm{red}_{t}^{-1}:=\redu^{-1}_{|SpT_{2n}(\varpi_t)}:SpT_{2n}(\varpi_t)&\rightarrow SST_{2n}(\varpi_l)\quad\\
\mathbf{a}&\mapsto\rm{red}_{t}^{-1}(\mathbf{a}):=T_0(a_1)T_1\cdots(a_t) T_t,\nonumber
\end{align}
where \begin{align}\label{inverseredu:nofactors}
&T_0=(1,2,\dots l_1),\\
&(a_i)T_i=
\begin{cases}
 ({a}_i)( a_i+1,\dots, a_i+l_{i+1}),
  &
 \mbox{if } a_i\in 2\mathbb{Z}\\
 ({a}_i)( a_i+1+1,a_i+1+2,\dots, a_i+1+l_{i+1}),
 &\mbox{if } a_i\notin 2\mathbb{Z}\\
\end{cases},& 1\le i\le t, \label{inverseredu:nofactors2}
\end{align}
and
$l_1,\dots, l_t$, $l_{t+1}$ are nonnegative even numbers defined by
\begin{align}\label{l:nofactors}
&l_1=
\begin{cases} min\{a_1-2, l-t\},& \mbox{if } a_1\in 2\mathbb{Z}\\
min\{ a_1-1, l-t\},&\mbox{if } a_1\notin 2\mathbb{Z},\\
\end{cases}
\end{align}
\begin{align}
&l_{i+1}=
\begin{cases}min\{a_{i+1}-a_i-2, l-t-\sum_{k=1}^i l_k\}, &\mbox{if }
a_i, a_{i+1}\in 2\mathbb{Z} \mbox{ or } a_i, a_{i+1}\notin 2\mathbb{Z}
\\
min\{{ a_{i+1}-a_i-1},~~ l-t-\sum_{k=1}^i l_k\},&\mbox{if }
a_i\in 2\mathbb{Z}, a_{i+1}\notin 2\mathbb{Z} \\
min\{a_{i+1}-a_i-3,  l-t-\sum_{k=1}^i l_k\},&\mbox{if }
 a_i\notin 2\mathbb{Z}, ~~~ a_{i+1}\in 2\mathbb{Z},
\end{cases},\quad \mbox{ for } 1\le i<t,\label{ells}
\end{align}
and $l_{t+1}=l-t-\sum_{i=1}^{t}l_i.$
\end{thm}

\begin{lem}\cite{azreco,azreduction}\label{lem:reductionsurj} Let $l\in [0,2n]$ and $t\in[n]\cap 2\Z$ such that  $2\le t\le min\{l,2n-l\}$ and
$l-t\in 2\Z$. Let
$a_1\notin 2\Z$,  and   $\mathbf{a} = (a_1, a_1+1,\dots , a_1+t-2,a_1+t-1)\in SpT_{2n}(\varpi_t)$. Then
  $a_1\ge t+1$,  and
\begin{align}\label{lem:inverseredx}
\rm{red}_{t}^{-1}:=\redu^{-1}_{|SpT_{2n}(\varpi_t)}:SpT_{2n}(\varpi_t)&\rightarrow SST_{2n}(\varpi_l)\quad\\
\mathbf{a}&\mapsto\rm{red}_{t}^{-1}(\mathbf{a})=T,\nonumber
\end{align}
where \begin{align}T=(1,2,\dots, l_1)\mathbf{a}(a_1+t-1+1,\dots, a_1+t-1+l_{t+1})\label{lem:T},
\end{align}
and $l_1$, $l_{t+1}$  are nonnegative even numbers defined by
\begin{align}\label{lem0}
l_1=min\{a_1-t-1,l-t\},\quad
l_{t+1}=l-t-l_1.
\end{align}
\end{lem}

\begin{thm}\cite{azreduction}\label{lem:reductionsurj2} Let $l\in [0,2n]$ and   $ t_i\neq 0\in 2\Z$, $i=1,\dots,k$, such that
$t=t_1+\cdots+t_k\in[n]$,
$$2\le t=\sum_{i=1}^kt_i\le min\{l,2n-l\}$$ and $l-t\in 2\Z$.
Let
\begin{align*}\mathbf{a} &=\bigoplus_{i=1}^k{\bf A_i}\in SpT_{2n}(\varpi_{t})
\end{align*}
where for $i=1,\dots,k$, ${a_i}\notin 2\Z$ and
\begin{align*}&{\bf A_i}=({a_i}, a_i+1,\dots , a_i+t_i-2,a_i+t_i-1).
\end{align*} Then $a_i\ge 2(t_1+\cdots+t_{i-1})+t_i+1$, $1\le i\le k$, and

\begin{align*}
\rm{red}_{t}^{-1}:=\redu^{-1}_{|SpT_{2n}(\varpi_t)}:SpT_{2n}(\varpi_t)&\rightarrow SST_{2n}(\varpi_l)\quad\\
\mathbf{a}&\mapsto\rm{red}_{t}^{-1}(\mathbf{a}):=T,\nonumber
\end{align*}
where \begin{align*}T&=(1,\dots, l_1){\bf {A_1}},( a_1+ 1+t_1-1,\dots,  a_1+ 1+t_1-1+l_2){\bf{A_2}}
( a_2+ 1+t_2-1,\dots,  a_2+ 1+t_2-1+l_3)\\
&\cdots{\bf A_{k-1}}( a_{k-1}+ 1+t_{k-1}-1,\dots,  a_{k-1}+ 1+t_{k-1}-1+l_{k})\\
&{\bf {A_k}}( {a_k}+ 1+t_k-1,\dots,  a_k+ 1+t_k-1+l_{t+1})
\end{align*}
with
$l_1,\dots, l_k$, $l_{t+1}$  nonnegative even numbers recursively defined for $i=1,\dots,k$, by
\begin{align}
&l_i=min\{~\displaystyle (a_j-2\sum_{m=1}^{j-1}t_m-t_{j}-1)-\sum_{m=1}^{i-1}l_m,~ i\le j\le k,~~l-t-\sum_{m=1}^{i-1}l_m\},\\
& \mbox{ and}\nonumber\\
&l_{t+1}=l-t-\sum_{i=1}^{k}l_i.
\end{align}
\end{thm}

\begin{thm}\cite{azreduction}\label{thm:xAy} Let $l\in [0,2n]$ and $t_1+p_0+p_1\in[n]$, $t_1\in 2\Z$ such that  $2\le t_1+p_0+p_1\le
min\{l,2n-l\}$ and $l-(t_1+p_0+p_1)\in 2\Z$. Let
$$\mathbf{a} = x_1x_2\cdots x_{p_0} {\bf A_1}y_1\cdots y_{p_1}\in SpT_{2n}(\varpi_{t_1+p_0+p_1})$$
 where $a_1\notin 2\Z$, ${\bf A_1}=(a_1, a_1+1,\dots , a_1+t_1-2,a_1+t_1-1)$ and
 $$\s(x_1x_2\cdots x_{p_0}y_1\cdots y_{p_1})\cap\mathbf{a}=\emptyset .$$ Then
  $a_1\ge t_1+2p_0+1$,  and
\begin{align}
\rm{red}_{t_1+p_0+p_1}^{-1}:=\redu^{-1}_{|SpT_{2n}(\varpi_{t_1+p_0+p_1})}:SpT_{2n}(\varpi_{t_1+p_0+p_1})&\rightarrow
SST_{2n}(\varpi_l)\nonumber\\
\mathbf{a}&\mapsto\rm{red}_{t_1+p_0+p_1}^{-1}(\mathbf{a})=T,\nonumber
\end{align}
where \begin{align}T=T_0(x_1)T_1\cdots  (x_{p_0})T_{p_0}\mathbf{A_1}Y_0(y_1)Y_1\cdots (y_{p_1})Y_{p_1}
\end{align}
is such that
\begin{enumerate}
\item  $\tilde l_1$ 
is the nonnegative even number defined by
\begin{align}
\tilde l_1=min\{a_1-t_1-2p_0-1,l-t_1-p_0-p_1\},\quad
\end{align}

\begin{align}\label{inverseredu:nofactors+++}
&T_0=(1,2,\dots l_1),\\
&({x}_j)T_j=
\begin{cases}
 ({x}_j)( x_j+1,\dots, x_j+l_{j+1}),
  &
 \mbox{if } x_j\in 2\mathbb{Z}\\
 ({x}_j)( x_j+1+1,x_j+1+2,\dots, x_j+1+l_{j+1}),
 &\mbox{if } x_j\notin 2\mathbb{Z}\\
\end{cases},& 1\le j\le p_0, \label{inverseredu:nofactors2++++}
\end{align}
with
$l_1,\dots, l_{p_0}$, $l_{p_0+1}$  nonnegative even numbers defined by
\begin{align}\label{l:nofactors++++}
&l_1=
\begin{cases} min\{x_1-2, \tilde l_1\},& \mbox{if } x_1\in 2\mathbb{Z}\\
min\{ x_1-1, \tilde l_1\},&\mbox{if } x_1\notin 2\mathbb{Z},\\
\end{cases}
\end{align}
\begin{align}
&l_{j+1}=
\begin{cases}min\{x_{j+1}-x_j-2, \tilde l_1-\sum_{m=1}^j l_m\}, &\mbox{if }
x_j, x_{j+1}\in 2\mathbb{Z} \mbox{ or } x_j, x_{j+1}\notin 2\mathbb{Z}
\\
min\{{ x_{j+1}-x_j-1},~~ \tilde l_1-\sum_{m=1}^j l_m\},&\mbox{if }
x_j\in 2\mathbb{Z}, x_{j+1}\notin 2\mathbb{Z} \\
min\{x_{j+1}-x_j-3,  \tilde l_1-\sum_{m=1}^j l_m\},&\mbox{if }
 x_j\notin 2\mathbb{Z}, ~~~ x_{j+1}\in 2\mathbb{Z},
\end{cases},\quad \mbox{ for } 1\le j<p_0,\label{ells+++}
\end{align}
and \begin{align}
 l_{p_0+1}=\tilde l_1-\sum_{m=1}^{p_0}l_m;\end{align}

\item put $y_0:=a_1+t_1-1\in 2\Z$, and define
\begin{align}\tilde l_{2}:=l-t_1-p_0-p_1- \sum_{m=1}^{p_0+1}l_m,
\end{align}

\begin{align}&Y_0=(y_0+1,y_0+2,\dots, y_0+ l_{ t_1+p_0+1}),\\
&({y}_j)Y_j=
\begin{cases}
 ({y}_j)( y_j+1,\dots, y_j+ l_{t_1+p_0+j+1}),
  &
 \mbox{if } y_j\in 2\mathbb{Z}\\
 ({y}_j)( y_j+1+1,y_j+1+2,\dots, y_j+1+ l_{t_1+p_0+j+1}),
 &\mbox{if } y_j\notin 2\mathbb{Z}\\
\end{cases},& 1\le j\le p_1, \label{inverseredu:nofactors2+++}
\end{align}
with
$ l_{t_1+p_0+j+1}$, $0\le j\le p_1$,
 nonnegative even numbers defined by
\begin{align}\label{l:nofactors+++}
& l_{t_1+p_0+1}=
\begin{cases} min\{y_1-y_0-2, \tilde l_2\},& \mbox{if } y_1\in 2\mathbb{Z}\\
min\{ y_1-y_0-1, \tilde l_2\},&\mbox{if } y_1\notin 2\mathbb{Z},\\
\end{cases}
\end{align}
\begin{align}
& l_{t_1+p_0+j+1}=
\begin{cases}min\{y_{j+1}-y_j-2, \tilde l_2-\sum_{m=1}^j  l_{t_1+p_0+m}\}, &\mbox{if }
y_j, y_{j+1}\in 2\mathbb{Z} \mbox{ or } y_, y_{j+1}\notin 2\mathbb{Z}
\\
min\{{ y_{j+1}-y_j-1},~~ \tilde l_2-\sum_{m=1}^j  l_{t_1+p_0+m}\},&\mbox{if }
y_j\in 2\mathbb{Z}, y_{j+1}\notin 2\mathbb{Z} \\
min\{y_{j+1}-y_j-3,  \tilde l_2-\sum_{m=1}^j  l_{t_1+p_0+m}\},&\mbox{if }
 y_j\notin 2\mathbb{Z}, ~~~ ,y_{j+1}\in 2\mathbb{Z},
\end{cases},\quad \mbox{ for } 1\le j<p_1,\label{ells++++}
\end{align}
and

\begin{align} l_{t_1+p_0+p_1+1}=\tilde l_2-\sum_{m=1}^{p_1}l_{t_1+p_0+m}.\end{align}
 \qquad\qquad\qquad\qquad\qquad\qquad\qquad\qquad\qquad\qquad\qquad\qquad\qquad\qquad\qquad\qquad\qquad\qquad\qquad
 \qquad$\Box$
\end{enumerate}
\end{thm}
A complete explicit formula generalizing the previous theorem is given in \cite{azreduction}.

Corollaries below are particular cases of the previous theorem with ${\bf A_1}=(b_1,\dots,b_{t_1})$ and they will be needed in the next section.
Next one considers   $p_0=1$, $p_1=0$ and $x_1=b_1-1\in 2\Z$.

\begin{cor}\label{cor:xA}Let $l\in [0,2n]$ and $t_1+1\in[n]$, $2\le t_1\in 2\Z$ such that  $ t_1+1\le min\{l,2n-l\}$ and $l-(t_1+1)\in 2\Z$.
Let
$$\mathbf{a} = (b_1-1) {\bf A_1}\in SpT_{2n}(\varpi_{t_1+1})$$
 where $b_1\notin 2\Z$, ${\bf A_1}=(b_1, b_1+1,\dots , b_1+t_1-2,b_1+t_1-1)$. 
 Then
  $b_1\ge t_1+2+1=t_1+3$, and
  \begin{align}\redu_{t_1+1}^{-1}(\mathbf{a})=T=(1,2,\dots,l_1)(b_1-1)\mathbf{A_1}((b_1+t_1-1)+1,\dots,(b_1+t_1-1)+l-(t_1+1)-l_1)\in
  SST_{2n}(\varpi_l)
  \end{align}
  where 
  $l_1=min\{b_1-(t_1+3), l-(t_1+1)\}$.
\end{cor}
\begin{proof}Let $p_0=1$, $p_1=0$ and $x_1=b_1-1\in 2\Z$. Since $b_1-1\in 2\Z$,  from Theorem \ref{thm:xAy} one has the nonnegative even
numbers

$$l_1=min\{(b_1-1)-2, \tilde l_1\}=\tilde l_1=min\{b_1-(t_1+3), l-(t_1+1)\}$$
and $l_2=\tilde l_1-l_1=0$.

Let $y_0=b_1+t_1-1$, since $p_1=0$,
$$l_{t_1+1+1}=\tilde l_2=l-t_1-1-l_1$$ and
$$Y_0=(y_0+1,y_0+2,\dots, y_0+ l_{ t_1+1+1})=((b_1+t_1-1)+1,\dots,(b_1+t_1-1)+l-t_1-1-l_1)\in SST_{2n}(l-t_1-1-l_1).$$
\end{proof}
Now one considers  $p_0=0$, $p_1=1$ and $y_1=b_{t_1}+1=b_1+t_1\notin 2\Z$.

\begin{cor}\label{cor:Ay}Let $l\in [0,2n]$ and $t_1+1\in[n]$, $t_1\in 2\Z$ such that  $2\le t_1+1\le min\{l,2n-l\}$ and $l-(t_1+1)\in 2\Z$. Let
$$\mathbf{a} =  {\bf A_1}(b_1+t_1)\in SpT_{2n}(\varpi_{t_1+1})$$
 where $b_1\notin 2\Z$, ${\bf A_1}=(b_1, b_1+1,\dots , b_1+t_1-2,b_1+t_1-1)$. 
 Then
  $b_1\ge t_1+1$, and
  \begin{align} \redu_{t_1+1}^{-1}(\mathbf{a})=T=(1,2,\dots,l_1)\mathbf{A_1}(b_1+t_1)(b_1+t_1+2,\dots,b_1+t_1+1+(l-t_1-1-l_1))\in SST_{2n}(\varpi_l)
  \end{align}
  where $l_1=min\{b_1-(t_1+1),l-(t_1+1)\}$.
\end{cor}

\medskip
Next one considers, $p_0=p_1=1$ and $x_1=b_1-1\in 2\Z$, $y_1=b_{t_1}+1=b_1+t_1\notin 2\Z$.

\begin{cor}\label{cor:xAy}Let $l\in [0,2n]$ and $t_1+2\in[n]$, $t_1\in 2\Z$ such that  $2\le t_1+2\le min\{l,2n-l\}$ and $l-(t_1+2)\in 2\Z$.
Let
$$\mathbf{a} = (b_1-1) {\bf A_1}(b_1+t_1)\in SpT_{2n}(\varpi_{t_1+2})$$
 where $b_1\notin 2\Z$, ${\bf A_1}=(b_1, b_1+1,\dots , b_1+t_1-2,b_1+t_1-1)$. 
 Then
  $b_1\ge t_1+2+1=t_1+3$, and
  \begin{align} \redu_{t_1+2}^{-1}(\mathbf{a})=T=(1,2,\dots,l_1)(b_1-1)\mathbf{A_1}(b_1+t_1)(b_1+t_1+2,\dots,b_1+t_1+1+(l-t_1-2-l_1))\in SST_{2n}(\varpi_l)
  \end{align}
  where $l_1=min\{b_1-(t_1+3), l-(t_1+2)\}$.
\end{cor}
\begin{proof}Let $p_0=1=p_1$ and $x_1=b_1-1\in 2\Z$ and $y_1=b_1+t_1\notin 2\Z$. Since $b_1-1\in 2\Z$,  from Theorem \ref{thm:xAy} one has the nonnegative even
numbers

$$l_1=min\{(b_1-1)-2, \tilde l_1\}=\tilde l_1=min\{b_1-(t_1+3), l-(t_1+2)\}$$
and $l_2=\tilde l_1-l_1=0\Rightarrow T_1=\emptyset$.

Let  $y_0=b_1+t_1-1\in 2\Z$
and
$$\tilde l_2=l-t_1-2-l_1-l_2=l-t_1-2-l_1$$.

 Then $y_1-y_0-1=0$ and  since $p_1=1$ and $y_1\notin 2\Z$,
$$l_{t_1+1+1}=min\{y_1-y_0-1, \tilde l_2\}=0\Rightarrow Y_0=\emptyset$$
Therefore,  $$l_{t_1+3}=\tilde l_2-l_{t_1+2}=\tilde l_2=l-t_1-2-l_1.$$

and
$$Y_1=(y_1+2,\dots, y_1+1+ l_{ t_1+3})=((b_1+t_1)+2,\dots,(b_1+t_1)+1+l-t_1-2-l_1)\in SST_{2n}(l-t_1-2-l_1).$$
\end{proof}

\begin{cor}\label{cor:xAystar}Let $l\in [0,2n]$ and $t_1+2\in[n]$, $t_1\in 2\Z$ such that  $2\le t_1+2\le min\{l,2n-l\}$ and $l-(t_1+2)\in 2\Z$.
Let 
$1\le u<b_1<b_1+t_1-1<b_1+t_1-1+v\le 2n$, $v\in\Z^+$,
 and
$$\mathbf{a} = (b_1-u) {\bf A_1}(b_1+t_1-1+v)\in SpT_{2n}(\varpi_{t_1+2})$$
 where $b_1\notin 2\Z$ and ${\bf A_1}=(b_1, b_1+1,\dots , b_1+t_1-2,b_1+t_1-1)$. 
 Then
  $b_1\ge t_1+2+1=t_1+3$, and
  \begin{align}& \redu_{t_1+2}^{-1}(\mathbf{a})=(1,2,\dots,l_1)\sqcup (b_1-u)\sqcup T_1\sqcup \mathbf{A_1}\sqcup\nonumber\\
   &\sqcup Y_0\sqcup (b_1+t_1-1+v)\sqcup(b_1+t_1+2,\dots,b_1+t_1+1+(l-t_1-2-l_1))\nonumber\\
  &\in SST_{2n}(\varpi_l)
  \end{align}
  where $l_1=min\{b_1-(t_1+3), l-(t_1+2)\}$.
\end{cor}

\section{Symplectic columns by reverse Schensted column bumping, the  action of the inverse reduction map and concatenation
}\label{sec:essential}

Fix $n\in\N$.
The purpose of this subsection is to analyse  when the action of the inverse reduction map
on a reverse bumped symplectic column $\b$ from a semi-standard tableau $S\in SST_{2n}(\mu)$ 
 returns a column $\redu^{-1}(\b)$ whose concatenation with the returned semi-standard tableau  $(  \longleftarrow S)$, via the reverse
 Schensted column bumping on $S$,  is a semi-standard tableau. This is   done  below for symplectic bumped intervals, not necessarily fixed by the parity involution,  in Subsection \ref{bumpintervalsymplecticinterval} and  Subsection \ref{bumpedinterval}, and for symplectic bumped columns disjoint of the action of the parity involution in Subsection  \ref{bumpedcolumndisjoint} that we recall from  \cite[Section 4.3]{azreco}.

Given the vector $(r_{1}<\dots<r_{m})$ with entries in $[1,u]$, its \emph{index} \emph{vector} $\delta_{(r_{1}<\dots<r_{m})}=(x_s)_{s\in
[1,u]}$  with respect to the interval $[1,u]$, is the vector in $\{0,1\}^u$ such that $x_s=1$ if $s\in\{r_1,\dots,r_m\}$ and $0$, otherwise.

Let $ m\in [0,n]$ and $S\in SST_{2n}(\mu)$    with $0\le m\le \ell(\mu)=t\le 2n$.
Let 
$r_{t_1}<\dots<r_{t_m}$, with $ \{t_1<\dots<t_m\}$ $\subseteq $ $[1,t]$, be $m$ row coordinates of   $S$ such that
$\mu-\delta_{(r_{t_1},\dots,r_{t_m})}$ is a partition. Let

$$(\b, \underset{ (r_{t_1},\dots,r_{t_m})}\longleftarrow S),$$
where
\begin{align*}\b=(a_{{t_1}}, \dots,a_{{t_m}})\in SpT_{2n}(\varpi_m),
\end{align*}
and  $$(\underset{(r_{t_1},\dots,r_{t_m})}\longleftarrow S)\in SST_{2n}(\mu-\delta_{(r_{t_1},\dots,r_{t_m})}).$$
For $ 1\le m\le t\le l\le 2n$,
such that
 $ m\le min\{l,2n-l\}$ and $l-m\in 2\Z$, we analyse  conditions under which one has
\begin{align}\rm{red}_{m}^{-1}(\mathbf{b})\bigcdot (\underset{(r_{t_1},\dots,r_{t_m})}\longleftarrow S) \in
SST_{2n}(\varpi_l+\mu-\delta_{(r_{t_1},\dots,r_{t_m})}).
\end{align}

Let $u=(u_1,\dots,u_t)\in \Z_+^t$ and $(v_1,\dots,v_{t'})\in \Z_+^{t'}$. We write $u\ge v$ whenever $u_i\ge v_i$ for $i=1,\dots,min\{t,t'\}$.
\begin{lem}\label{lem:symp} If  ${\bf a}=(a_1,\dots,a_t)\in SpT_{2n}(\varpi_t)$ and ${\bf b}=(b_1,\dots,b_{t'})\in SST_{2n}(\varpi_{t'})$, for
some $t'\le t$, are such that ${\bf b} \ge {\bf a}$, that is, $b_i\ge a_i$, $1\le i\le t'$,  then ${\bf b}\in SpT_{2n}(\varpi_{t'})$.
\end{lem}
\begin{proof}If ${\bf a}=(a_1,\dots,a_t)\in SpT_{2n}(\varpi_t)$ and ${\bf b} \ge {\bf a}$, then $b_k\ge a_k\ge 2k-1$, for $k=1,\dots,t'$.
Hence, ${\bf b}\in SpT_{2n}(\varpi_{t'})$.
\end{proof}
\begin{lem}\label{keybasic} Let $m\in \N$ and  $(a_1< \dots<a_{m-1}< a_{m})\in\Z^m_{\ge 1}$. Then
\begin{itemize}
\item[(a)] $(1,2,\dots,m-1,m)\le (a_1, \dots,a_{m-1}, a_{m})$
and
\item[(b)]$ (1,2,\dots,m-1,l)\bigcdot (a_1, \dots,a_{m-1}, a_{m})\in SST(\varpi_l+\varpi_m)$ is a two-column semi-standard tableau of length
    $l$ for any $l\ge m$.
\item[(c)] for $p\ge a_m-1$,
$$ (1,2,\dots,p, p+1)\bigcdot (a_1, \dots,a_{m-1}, a_{m})\in SST(\varpi_{p+1}+\varpi_m),  $$
is a two column semi-standard tableau of length $p+1\ge m$.
\end{itemize}
\end{lem}
\begin{proof} $(a)$ One has $1\le a_1$. Let $k\ge 1$. Assuming by induction that $(1,2,\dots,k)\le  (a_1, \dots,a_{k})$ it follows that
$k+1\le a_k+1\le a_{k+1}$ and $ (1,2,\dots,k+1)\le (a_1, \dots,a_{k+1})$.

$(b)$ Follows from $(a)$.

$(c)$ It follows from $(b)$ since $p\ge a_m-1\ge a_{m-1}\ge m-1\Rightarrow p+1\ge m$.
\end{proof}


Recall the reverse Schensted column insertion in Section \ref{subsec:insert}, 
and explicit formula for the inverse
reduction map  in   Lemma \ref{lem:reductionsurj}, Theorem  \ref{thm:nofactors} and Theorem \ref{thm:xAy}, with their corollaries, and more generally in
\cite{azreduction}.
\subsection{Symplectic bumped columns disjoint of the action of the parity involution }\label{bumpedcolumndisjoint}
{ The following theorem essentially appears in \cite[Section 4.3]{azreco}  and is stated here with a bit more generality. It  considers the case where the symplectic bumped  column $\b$ satisfies $\s(\b)\cap\b=\emptyset$.}

\begin{thm}\cite{azreco} \label{thm:essentialss}Let $ z\in [n]$ and $S\in SST_{2n}(\mu)$    
 with first column
 $\mathbf{a} = (a_1, \dots , a_t)$ $\in$ $ SST_{2n}(\varpi_t)$. 
Let  $r_{t_1}<\dots<r_{t_z}$ with $ \{t_1<\dots<t_z\}\subseteq [1,t]$ be $z$ row coordinates of   $S$ such that
$\mu-\delta_{(r_{t_1},\dots,r_{t_z})}$ is a partition. Let

$$(\b, \underset{ (r_{t_1},\dots,r_{t_z})}\longleftarrow S),$$
where
\begin{align*}\b=(a_{{t_1}}, \dots,a_{{t_z}})\in SpT_{2n}(\varpi_z)  \mbox{ is such that $\s(\b)\cap \bf{b}=\emptyset$},
\end{align*}
and  $$(\underset{(r_{t_1},\dots,r_{t_z})}\longleftarrow S)\in SST_{2n}(\mu-\delta_{(r_{t_1},\dots,r_{t_z})})$$
has first column $S^1=(a_1, \dots,a_{t_1-1},  x_{t_1},a_{t_1+1},\dots ,a_{t_z-1},x_{t_z},a_{t_z+1},\dots, a_{t})\in SST_{2n}(\varpi_t)$.
Then
\begin{enumerate}
\item 
$\b=(a_{t_1}, \dots,a_{t_z})\le (x_{t_1},\dots, x_{t_z}).$
\item For $\oo z\le n$, $1\le z\le t\le l\le 2n$ 
such that
 $ z\le min\{l,2n-l\}$ and $l-z\in 2\Z$,
 \end{enumerate}


 \begin{align}\redu_z^{-1}(a_{t_1}, \dots,a_{t_z})=\rm{red}_{z}^{-1}(\mathbf{b}):=T_0(a_{t_1})T_1\cdots(a_{t_z}) T_z\in
 SST_{2n}(\varpi_l),\nonumber
\end{align}
where \begin{align}
&T_0=(1,2,\dots l_1),\\
&(a_{t_i})T_i=
\begin{cases}
 ({a}_{t_i})( a_{t_i}+1,\dots, a_{t_i}+l_{i+1}),
  &
 \mbox{if } a_{t_i}\in 2\mathbb{Z}\\
 ({a}_{t_i})( a_{t_i}+1+1,a_{t_i}+1+2,\dots, a_{t_i}+1+l_{i+1}),
 &\mbox{if } a_{t_i}\notin 2\mathbb{Z}\\
\end{cases},& 1\le i\le z,
\end{align}
and
$l_1,\dots, l_z$, $l_{z+1}$ are nonnegative even numbers defined by
\begin{align}
&l_1=
\begin{cases} min\{a_{t_1}-2, l-z\},& \mbox{if } a_{t_1}\in 2\mathbb{Z}\\
min\{ a_{t_1}-1, l-z\},&\mbox{if } a_{t_1}\notin 2\mathbb{Z},\\
\end{cases}
\end{align}
\begin{align}
&l_{i+1}=
\begin{cases}min\{a_{t_{i+1}}-a_{t_i}-2, l-z-\sum_{k=1}^i l_k\}, &\mbox{if }
a_{t_i}, a_{t_{i+1}}\in 2\mathbb{Z} \mbox{ or } a_{t_i}, a_{t_{i+1}}\notin 2\mathbb{Z}
\\
min\{ a_{t_{i+1}}-a_{t_i}-1,~~ l-z-\sum_{k=1}^i l_k\},&\mbox{if }
a_{t_i}\in 2\mathbb{Z}, a_{t_{i+1}}\notin 2\mathbb{Z} \\
min\{a_{t_{i+1}}-a_{t_i}-3,  l-z-\sum_{k=1}^i l_k\},&\mbox{if }
 a_{t_i}\notin 2\mathbb{Z}, ~~~ a_{t_{i+1}}\in 2\mathbb{Z},
\end{cases},\quad \mbox{ for } 1\le i<z,\label{ellsagain}
\end{align}
and $\displaystyle l_{z+1}=l-z-\sum_{i=1}^{z}l_i.$

Moreover,
\begin{align*}&T_0(a_{t_1})\bigcdot (a_1, \dots,a_{t_1-1},  x_{t_1})\in SST_{2n}(\varpi_{l_1+1}+\varpi_{t_1})\\
&T_1(a_{t_2})\bigcdot(a_{t_1+1},\dots ,a_{t_2-1},x_{t_2})
\in SST_{2n}(\varpi_{l_2+1}+\varpi_{t_2-t_1})\\
&\vdots\\
&T_{z-1}(a_{t_z})\bigcdot(a_{t_{z-1}+1},\dots ,a_{t_z-1},x_{t_z})
\in SST_{2n}(\varpi_{l_{z}+1}+\varpi_{t_z-t_{z-1}})\\
&T_z\bigcdot (a_{t_z+1},\dots, a_{t})\in SST_{2n}(\varpi_{l_{z+1}}+\varpi_{t-t_z})\\
\end{align*}
and
\begin{align}T_0(a_{t_1})T_1\cdots(a_{t_z}) T_z\bigcdot (a_1, \dots,a_{t_1-1},  x_{t_1},a_{t_1+1},\dots ,a_{t_z-1},x_{t_z},a_{t_z+1},\dots,
a_{t})\in SST_{2n}(\varpi_l+\varpi_t).
\end{align}
\begin{enumerate}
\item [(3)]$\rm{red}_{z}^{-1}(\mathbf{b})\bigcdot (\underset{(r_{t_1},\dots,r_{t_z})}\longleftarrow S)$ $\in
    SST_{2n}(\varpi_l+\mu-\delta_{(r_{t_1},\dots,r_{t_z})})$.
    \end{enumerate}
\end{thm}

\begin{obs}  In the previous theorem, we have considered $(\underset{(r_{t_1},\dots,r_{t_z})}\longleftarrow S)$ with the same length as $S$.
 When the length of $S$ is shortened by the reverse  bumping,  the proof is similar.
\end{obs}

\subsection{Bumped  columns from symplectic intervals closed for the parity  involution}\label{bumpintervalsymplecticinterval}
 We consider (symplectic) bumped columns from
a symplectic tableau where the first column is a symplectic interval $\a$ such that $\s(\a)=\a$.

\subsubsection{The symplectic bumped column is disjoint of the action of the parity involution}
If the symplectic bumped column $\b$ is such that $\b\cap \s(\b)=\emptyset$ then  we are in the case of Subsection \ref{bumpedcolumndisjoint}.

 \subsubsection{The bumped column is an interval}We first consider  a symplectic interval $\a$ such that $\s({\a})={\a}$ from where the  interval  $\b$  not necessarily fixed by the parity involution is  bumped.
The statement appears in \cite{azreco}. We detail a bit more its proof. We have to consider two cases: the bumped interval $\b$ from $\a$ starts   with  an odd
number or with an even number.

\begin{prop} \cite{azreco}\label{reducbumpinterval1} 
Let $k+t\in[n]\cap 2\Z$, for some $ k\ge 0$ and $t\ge 1$,
$a_1\notin 2\Z$,  and   $\mathbf{a} = (a_1, a_1+1,\dots ,a_1+{k-1},a_1+k, a_1+k+1,\dots,a_1+k+t-1)\in SpT_{2n}(\varpi_{k+t})$ with
(necessarily) $a_1\ge k+t+1$.
Let
$$(\b, \underset{ (k+1,\dots,k+t)}\longleftarrow \a),$$
where
$$\b=(a_1+k, a_1+k+1,\dots,a_1+k+t-1)  \mbox{ and  }(\underset{(k+1,\dots,k+t)}\longleftarrow \a)=(a_1, a_1+1,\dots ,a_1+{k-1})\in
SpT_{2n}(\varpi_{k}).$$
Then
\begin{enumerate}
\item $\b=(a_1+k, a_1+k+1,\dots,a_1+k+t-1)\in SpT_{2n}(\varpi_{t})$.
\item for $k+t\le l\le 2n$, 
such that
 $ t\le min\{l,2n-l\}$ and $l-t\in 2\Z$,

\medskip
$(a)$ if $k,t\in[n]\cap 2\Z\Leftrightarrow$ $ a_1+k\notin 2\Z$, $a_1+k+t-1\in 2\Z$, it holds $\s(\b)=\b$ and

\begin{align} \label{redub}
{\redu}_{t}^{-1}(\mathbf{b})=(1,2,\dots, l_1){\b}((a_1+k-1+t)+1,\dots, (a_1+k-1+t)+(l-t-l_1))\in SST_{2n}(\varpi_l)
\end{align}
where 
$l_1$, $l_{t+1}$  are nonnegative even numbers defined by
\begin{align*}
k\le l_1=min\{a_1+k-t-1,l-t\},\quad
l_{t+1}=l-t-l_1,
\end{align*}
and the  concatenation of $\redu_{t}^{-1}(\b) \in SST_{2n}(\varpi_l)$ and $(\underset{(k+1,\dots,k+t)}\longleftarrow \a)\in
SpT_{2n}(\varpi_{k})$ 
is a semi-standard tableau in $ SST_{2n}(\varpi_l+\varpi_k)$, with $l\ge k+t$. That is,
 \begin{align}
 &\redu_{t}^{-1}(\b)\bigcdot (\underset{(k+1,\dots,k+t)}\longleftarrow \a)=\nonumber\\
 &={\redu}_{t}^{-1}(\mathbf{b})\bigcdot (a_1, a_1+1,\dots ,a_1+{k-1})\nonumber\\
 &=(1,2,\dots, l_1,a_1+k, a_1+k+1,\dots,a_1+k+t-1,a_1+k+t-1+1,\dots, a_1+k+t-1+(l-t-l_1))\bigcdot\nonumber\\
 &\bigcdot(a_1, a_1+1,\dots ,a_1+{k-1})\in SST_{2n}(\varpi_l+\varpi_k),\; l\ge k+t,\; l_1\ge k.
 \end{align}
 In other words the row coordinates of the bumped entries in $\a$ and defining $\b$,  are pushed down or fixed in $\redu_{t}^{-1}(\b)$.
with $l\ge k+t$.

{ In particular, when $k=0$ and  $2\le t\in 2\Z$, then $\b=\a$ and $(\underset{(1,\dots,t)}\longleftarrow \a)=( )\in SpT_{2n}()$.
Moreover,
\begin{align*}&{\redu}_{t}^{-1}(\mathbf{b})=\\
 &=(1,2,\dots, l_1,a_1, a_1+1,\dots,a_1+t-1,a_1+t-1+1,\dots, a_1+t-1+(l-l_1-t))\in SST_{2n}(\varpi_l),
 \end{align*}
 \noindent with  $l_1=min\{a_1-t-1,l-t\}$, $a_1\ge t+1.$ }

\medskip
$(b)$ if $k,t\in [n]$, and $k,t \notin 2\Z$ $\Leftrightarrow a_1+k\in 2\Z,          \, a_1+k+t-1\in 2\Z$, it holds

 $$\b=(a_1+k)\sqcup (a_1+k+1,\dots,a_1+k+t-1)\in SpT_{2n}(\varpi_{t}), \;a_1+k+1\notin 2\Z$$
 $$=(a_1+k)\sqcup{\overline \b},\;\; \overline \b\in  SpT_{2n}(\varpi_{t-1}), \; \s(\overline \b)=\overline \b.$$

 If $\overline\b=\emptyset\Leftrightarrow t=1$ and $a_1\ge k+2$,
 \begin{align}\label{b=empty}
{\redu}_{t}^{-1}(\mathbf{b})=(1,2,\dots, l_1)\sqcup (a_1+k)\sqcup (a_1+k+1,\dots, a_1+k+1+(l-t-l_1))\in SST_{2n}(\varpi_l),
\end{align}
where
\begin{align}
&k\le l_1=
min\{a_1+k-2, l-1\}, 
\end{align}
and \begin{align}\label{notemptyt=1}&{\redu}_{t}^{-1}(\mathbf{b})\bigcdot(\underset{(k+1)}\longleftarrow \a)=\nonumber\\
 &{\redu}_{t}^{-1}(\mathbf{b})\bigcdot (a_1, a_1+1,\dots ,a_1+{k-1})=\nonumber\\
 &=(1,2,\dots, l_1,a_1+k,a_1+k+1,\dots, a_1+k+1+l-1-l_1)\bigcdot (a_1, a_1+1,\dots ,a_1+{k-1})\\
 &\in SST_{2n}(\varpi_l+\varpi_k),\; l\ge k+1,\, l_1\ge k.\nonumber
\end{align}

\medskip
Otherwise, $\overline\b\neq\emptyset\Leftrightarrow 2\le t-1\in 2\Z$
\begin{align}\label{notempty} &\redu_{t}^{-1}(\mathbf{b})=\nonumber\\
&=(1,2,\dots, l_1)\sqcup(a_1+k){\overline\b}\sqcup((a_1+k-1+t)+1,\dots,
a_1+k-1+t+(l-(t-1)-1-l_1))\in SST_{2n}(\varpi_l),\nonumber\\
&l_1\ge k,
\end{align}
where
\begin{align}
& k\le l_1=min\{a_1+k+1-(t-1+3), l-(t-1+1)\} 
\end{align}
and \begin{align}&{\redu}_{t}^{-1}(\mathbf{b})\bigcdot(\underset{(k+1,\dots,k+t)}\longleftarrow \a)=\nonumber\\
&={\redu}_{t}^{-1}(\mathbf{b})\bigcdot (a_1, a_1+1,\dots ,a_1+{k-1})\in  SST_{2n}(\varpi_l+\varpi_k),\; l\ge k+t,\, l_1\ge k.
\end{align}
\end{enumerate}
\end{prop}
\begin{proof}
$(1)$ Follows from 
Lemma  \ref{lem:symp}.

$(2)$ Note $l-t\ge k+t-t\ge k$ and $a_1+k-t-1\ge k+t+1+k-t-1=2k$.

$(a)$. From Lemma \ref{lem:reductionsurj}, ${\redu}_{t}^{-1}(\mathbf{b})$ is given by \eqref{redub}. Therefore $l_1\ge k$. This means that the
row coordinates of the entries of the bumped column $\b$ as entries of $\a$ have been pushed down or kept in $\redu_{t}^{-1}(\b)$,
equivalently  $(1,2,\dots, k)\le (a_1, a_1+1,\dots ,a_1+{k-1})$.

Since $\ell(\redu_{t}^{-1}(\b))=l$, the concatenation of $\redu_{t}^{-1}(\b)$ and $(a_1, a_1+1,\dots ,a_1+{k-1})$ is a semi-standard tableau
in $ SST_{2n}(\varpi_l+\varpi_k)$, with $l\ge k+t$.

$(b)$ If $\overline\b=\emptyset\Leftrightarrow t=1$, \eqref{b=empty} is given by Theorem \ref{thm:nofactors}.

If $\overline\b\neq\emptyset\Leftrightarrow 2\le t-1\in 2\Z$,  \eqref{notempty} is given by
  Theorem \ref{thm:xAy} in the form $(a_1+k)\b$  with $t_1=t-1$ and { $a_1+k\ge 2k+t+1$}. Then, from Corollary \ref{cor:xA}, one has
  with $b_1=a_1+k+1\notin 2\Z$

  $$\mathbf{b} = (b_1-1) {\overline\b}\in SpT_{2n}(\varpi_{t_1+1})$$
  $\overline\b=(b_1, b_1+1,\dots , b_1+t_1-2,b_1+t_1-1)$. 
 Then
  $b_1=a_1+k+1\ge 2k+t+2\ge 2k+ t_1+2+1=2k+t_1+3$, and
  \begin{align}\redu_{t_1+1}^{-1}(\mathbf{b})=(1,2,\dots,l_1)(b_1-1)\overline\b((b_1+t_1-1)+1,\dots,(b_1+t_1-1)+l-t_1-1-l_1)\in
  SpT_{2n}(\varpi_{l})
  \end{align}
  where $k\le l_1=min\{b_1-3,\tilde l_1\}=\tilde l_1=min\{b_1-(t_1+3), l-(t_1+1)\}$. Note

  $$ l_1=min\{(a_1+k+1)-3, \tilde l_1\}=\tilde l_1=min\{a_1+k+1-(t-1+3), l-(t-1+1)\}$$

and   $a_1+k+1-(t-1+3)\ge 2k+(t-1)+3-(t-1+3)=2k$,  and  $l-(t-1+1)\ge k+t-t\ge k$.

Hence \begin{align} &\redu_{t}^{-1}(\mathbf{b})=\nonumber\\
&=(1,2,\dots, l_1)\sqcup(a_1+k){\overline\b}\sqcup((a_1+k-1+t)+1,\dots,
(a_1+k-1+t)+(l-(t-1)-1-l_1))\in SST_{2n}(\varpi_l),\nonumber\\
&l_1\ge k.\nonumber
\end{align}
\end{proof}

The following is a generalization of the previous proposition for symplectic tableaux where the first column  $\a$ is is symplectic interval such that $\s({\a})={\a}$ and the bumped symplectic column $\b$ is an interval
  at the  bottom  of $\a$.

\begin{thm}\label{thm:reducbumpinterval1} 

Let $k+t\in[n]\cap 2\Z$, for some $ k\ge 0$ and $t\ge 1$,
$a_1\notin 2\Z$,  and  $S\in SpT_{2n}(\mu)$ a symplectic tableau with first column the symplectic interval  $\mathbf{a} = (a_1, a_1+1,\dots
,a_1+{k-1},a_1+k, a_1+k+1,\dots,a_1+k+t-1)$ $\in SpT_{2n}(\varpi_{k+t})$ such that $\s(\a)=\a$ and   (necessarily) $a_1\ge k+t+1$.
Let $r_{k+1}<\dots<r_{k+t}$ be $t$ row coordinates of   $S$ such that
$\mu-\delta_{(r_{k+1},\dots,r_{k+t})}$ is a partition, and
$$(\b, \underset{ (r_{k+1},\dots,r_{k+t})}\longleftarrow S),$$
where
$$\b=(a_1+k, a_1+k+1,\dots,a_1+k+t-1)\in SpT_{2n}(\varpi_t)$$   and  $$(\underset{(r_{k+1},\dots,r_{k+t})}\longleftarrow S)\in
SpT_{2n}(\mu-\delta_{(r_{k+1},\dots,r_{k+t})})$$
has first column $S^1=(a_1, a_1+1,\dots ,a_1+{k-1},x_k,x_{k+1}\dots, x_{k+t'-1})\in SpT_{2n}(\varpi_{k+t'}),$ for some $0\le t'\le t$.
Then
\begin{enumerate}
\item $\b=(a_1+k, a_1+k+1,\dots,a_1+k+t-1)\in SpT_{2n}(\varpi_t)$ and $S^1\in SpT_{2n}(\varpi_{k+t'})$ are such that

 $$(x_k,x_{k+1}\dots, x_{k+t'-1})\ge \b=(a_1+k, a_1+k+1,\dots,a_1+k+t'-1,\dots, a_1+k+t-1),\; 1\le  t'\le t,$$
\item For $k+t\le l\le 2n$, 
such that
 $ t\le min\{l,2n-l\}$ and $l-t\in 2\Z$:

\medskip
$(a)$ if $k,t\in[n]\cap 2\Z\Leftrightarrow$ $ a_1+k\notin 2\Z$, $a_1+k+t-1\in 2\Z$, it holds $\s(\b)=\b$ and

\begin{align} \label{redubthm}
{\redu}_{t}^{-1}(\mathbf{b})=(1,2,\dots, l_1){\b}(a_1+k-1+t+1,\dots, a_1+k-1+t+(l-t-l_1))\in SST_{2n}(\varpi_l),
\end{align}
where 
$l_1$, $l_{t+1}$  are nonnegative even numbers defined by
\begin{align*}
k\le l_1=min\{a_1+k-t-1,l-t\},\quad
l_{t+1}=l-t-l_1,
\end{align*}

and
 \begin{align}&{\redu}_{t}^{-1}(\mathbf{b})\bigcdot S^1=\nonumber\\
 &={\redu}_{t}^{-1}(\mathbf{b})\bigcdot (a_1, a_1+1,\dots ,a_1+{k-1},x_k,x_{k+1}\dots, x_{k+t'-1})=\nonumber\\
 &=(1,2,\dots, l_1,a_1+k, a_1+k+1,\dots,a_1+k+t-1,a_1+k+t-1+1,\dots, a_1+k+t-1+(l-t-l_1))\bigcdot\nonumber\\
 &\bigcdot(a_1, a_1+1,\dots ,a_1+{k-1},x_k,x_{k+1}\dots, x_{k+t'-1})\in SST_{2n}(\varpi_l+\varpi_{k+t'}),\nonumber\\
 & \mbox{ with $ l_1\ge k,\; l_1+t\ge k+t'$ and $l\ge l_1+t\ge k+t\ge k+t'$}.
 \end{align}
 That is, the  concatenation of $\redu_{t}^{-1}(\b)$ and $S^1=(a_1, a_1+1,\dots ,a_1+{k-1},x_k,x_{k+1}\dots, x_{k+t'-1})$ is a semi-standard
 tableau in $ SST_{2n}(\varpi_l+\varpi_{k+t'})$, with $l\ge k+t\ge k+t'$. Moreover,

 \begin{align}
 &{\redu}_{t}^{-1}(\mathbf{b})\bigcdot (\underset{ (r_{k+1},\dots,r_{k+t})}\longleftarrow S)\in
 SST_{2n}(\varpi_l+\mu-\delta_{(r_{k+1},\dots,r_{k+t})})
 \end{align}

{ In particular, when $k=0$ and $2\le t\in 2\Z$, $\b=\a$, and $(\underset{(1,\dots,t)}\longleftarrow S)\in SpT_{2n}(\mu-\varpi_t )$.
Moreover,
\begin{align*}&{\redu}_{t}^{-1}(\mathbf{b})=\\
 &=(1,2,\dots, l_1,a_1, a_1+1,\dots,a_1+t-1,a_1+t-1+1,\dots, a_1+t-1+l_{t+1})\in SST_{2n}(\varpi_l),\\
 &\; l=l_1+l_{t+1}+t ,
 \end{align*}
 \noindent with  $l_1=min\{a_1-t-1,l-t\}$, $a_1\ge t+1,$ and \begin{align}
 &{\redu}_{t}^{-1}(\mathbf{b})\bigcdot (\underset{ (1,\dots,t)}\longleftarrow S)\in SST_{2n}(\varpi_l+\mu-\delta_{(1,\dots,t)}).
 \end{align}}

\medskip
$(b)$ if $k,t\in [n]$, and $k,t \notin 2\Z$ $\Leftrightarrow a_1+k\in 2\Z,          \, a_1+k+t-1\in 2\Z$, it holds

 $$\b=(a_1+k)\sqcup (a_1+k+1,\dots,a_1+k+t-1)\in SpT_{2n}(\varpi_{t}), \;a_1+k+1\notin 2\Z$$
 $$=(a_1+k).{\overline \b},\;\; \overline \b\in  SpT_{2n}(\varpi_{t-1}), \, \s(\overline \b)=\overline \b.$$

 If $\overline\b=\emptyset\Leftrightarrow t=1$
 \begin{align}\label{b=emptythm}
{\redu}_{t}^{-1}(\mathbf{b})=(1,2,\dots, l_1)(a_1+k)(a_1+k+1,\dots, a_1+k+1+(l-1-l_1))\in SST_{2n}(\varpi_l),
\end{align}
where
\begin{align}
&k\le l_1=
min\{a_1+k-2, l-1\},
\end{align}

and \begin{align}
&{\redu}_{t}^{-1}(\mathbf{b})\bigcdot S^1=\nonumber\\
 &={\redu}_{t}^{-1}(\mathbf{b})\bigcdot (a_1, a_1+1,\dots ,a_1+{k-1},x_k)=\nonumber\\
 &=(1,2,\dots, l_1,a_1+k,a_1+k+1,\dots, a_1+k+1+l-1-l_1)\bigcdot (a_1, a_1+1,\dots ,a_1+{k-1},x_k)\\
 &\in SST_{2n}(\varpi_l+\varpi_k),\; l\ge k+1,\, l_1\ge k.\nonumber
\end{align}

Otherwise, $\overline\b\neq\emptyset\Leftrightarrow 2\le t-1\in 2\Z$,

\begin{align}& \redu_{t}^{-1}(\mathbf{b})=(1,2,\dots, l_1)(a_1+k){\overline\b}(a_1+k-1+t+1,\dots, a_1+k-1+t+(l-t-l_1))\in SST_{2n}(\varpi_l),\nonumber\\
&\;
l_1\ge k
\label{notemptythm}
\end{align}
where
\begin{align}
& k\le l_1=min\{a_1+k+1-(t-1+3), l-(t-1+1)\} 
\end{align}

and \begin{align*}&{\redu}_{t}^{-1}(\mathbf{b})\bigcdot S^1=\nonumber\\
&{\redu}_{t}^{-1}(\mathbf{b})\bigcdot (a_1, a_1+1,\dots ,a_1+{k-1},x_k,x_{k+1}\dots, x_{k+t'-1})\in  SST_{2n}(\varpi_l+\varpi_{k+t'}),\\
& l\ge k+t\ge k+t',\, l_1\ge k.
\end{align*}
\end{enumerate}

That is, the  concatenation of $\redu_{t}^{-1}(\b)$ and $S^1=(a_1, a_1+1,\dots ,a_1+{k-1},x_k,x_{k+1}\dots, x_{k+t'-1})$ is a semi-standard
tableau in $ SST_{2n}(\varpi_l+\varpi_{k+t'})$, with $l\ge k+t\ge k+t'$. Moreover,

 \begin{align}
 &{\redu}_{t}^{-1}(\mathbf{b})\bigcdot (\underset{ (r_{k+1},\dots,r_{k+t})}\longleftarrow S)\in
 SST_{2n}(\varpi_l+\mu-\delta_{(r_{k+1},\dots,r_{k+t})}).
 \end{align}
 $\qquad\qquad\qquad \qquad\qquad\qquad \qquad\qquad\qquad \qquad\qquad\qquad \qquad\qquad\qquad \qquad\qquad\qquad \qquad\qquad\qquad \Box$
\end{thm}

\medskip

Next theorem considers a symplectic tableau where the first column is a symplectic interval $\a$ such that $\s(\a)=\a$  and the bumped symplectic interval $\b$ from $\a$ is in general position. It follows from Lemma \ref{lem:reductionsurj}, Corollary \ref{cor:xA}, Corollary \ref{cor:Ay} and Corollary \ref{cor:xAy}.

\begin{thm}\label{thm:reducbumpinterval1m+} 

Let $k+t+s\in[n]\cap 2\Z$, for some $ k,s\ge 0$ and $t\ge 1$,
$a_1\notin 2\Z$. Let $S\in SpT_{2n}(\mu)$ be a symplectic tableau with first  column  the symplectic interval $\a$ in $ SpT_{2n}(\varpi_{k+t+s})$ such that $\s(\a)=\a$,  $$\mathbf{a}
= (a_1, a_1+1,\dots ,a_1+{k-1})\sqcup (a_1+k, a_1+k+1,\dots,a_1+k+t-1)\sqcup (a_1+k+t,a_1+k+t+1\dots,a_1+k+t+s-1)$$ where (necessarily)
$a_1\ge k+t+s+1$.
Let $r_{k+1}<\dots<r_{k+t}$ be $t$ row coordinates of   $S$ such that
$\mu-\delta_{(r_{k+1},\dots,r_{k+t})}$ is a partition, and
$$(\b, \underset{ (r_{k+1},\dots,r_{k+t})}\longleftarrow S),$$
with
$$\b=(a_1+k, a_1+k+1,\dots,a_1+k+t-1)\in SpT_{2n}(\varpi_t)$$   and  $$(\underset{(r_{k+1},\dots,r_{k+t})}\longleftarrow S)\in
SpT_{2n}(\mu-\delta_{(r_{k+1},\dots,r_{k+t})})$$
has first column $S^1\in SpT_{2n}(\varpi_{k+t+s}).$
Then
\begin{enumerate}
\item 
$S^1=\a\in SpT_{2n}(\varpi_{k+t+s})$.

\item For $k+t+s\le l\le 2n$, 
such that
 $ t\le min\{l,2n-l\}$ and $l-t\in 2\Z$:

\medskip
$(a)$ If $k,t,s\in[n]\cap 2\Z\Leftrightarrow$ $ a_1+k\notin 2\Z$, $ a_1+k+t\notin 2\Z$, 
it holds

\begin{align} \label{redubthm+}
{\redu}_{t}^{-1}(\mathbf{b})=(1,2,\dots, l_1){\b}(a_1+k-1+t+1,\dots, a_1+k-1+t+(l-t-l_1))\in SST_{2n}(\varpi_l),
\end{align}
\noindent where 
and $l_1$, $l_{t+1}$  are nonnegative even numbers defined by
\begin{align*}
k+s\le l_1=min\{a_1+k-t-1,l-t\},\quad
l_{t+1}=l-t-l_1,
\end{align*}
and
 \begin{align}&{\redu}_{t}^{-1}(\mathbf{b})\bigcdot S^1={\redu}_{t}^{-1}(\mathbf{b})\bigcdot\a \nonumber\\
 &=(1,2,\dots, l_1)\sqcup (a_1+k, a_1+k+1,\dots,a_1+k+t-1)\sqcup \nonumber\\
 &\sqcup(a_1+k+t-1+1,\dots, a_1+k+t-1+(l-t-l_1))\bigcdot\nonumber\\
 &\bigcdot(a_1, a_1+1,\dots ,a_1+{k-1})\sqcup (a_1+k, a_1+k+1,\dots,a_1+k+t-1)\sqcup (a_1+k+t,\dots, a_1+k+t+s-1)\nonumber\\
 &\in SST_{2n}(\varpi_l+\varpi_{k+t+s}),\nonumber\\
 & \mbox{ with $ l_1\ge k+s,\; l_1+t\ge k+s+t$ and $l=l_1+t+l_{t+1}\ge k+s+t+l_{t+1}\ge k+t+s$}.
 \end{align}
 That is, the  concatenation of $\redu_{t}^{-1}(\b)$ and $S^1$ 
 is a semi-standard tableau in $ SST_{2n}(\varpi_l+\varpi_{k+t+s})$, with $l\ge k+t+s$. Moreover,

 \begin{align}
 &{\redu}_{t}^{-1}(\mathbf{b})\bigcdot (\underset{ (r_{k+1},\dots,r_{k+t})}\longleftarrow S)\in
 SST_{2n}(\varpi_l+\mu-\delta_{(r_{k+1},\dots,r_{k+t})})
 \end{align}

\medskip
$(b)$  if $ k,t, s\in [n]$, and $s\in 2\Z$, $k,t \notin 2\Z$ $\Leftrightarrow a_1+k\in 2\Z, \, a_1+k+t-1\in 2\Z$,  it holds

 $$\b=(a_1+k)\sqcup (a_1+k+1,\dots,a_1+k+t-1)\in SpT_{2n}(\varpi_{t}), \;a_1+k+1\notin 2\Z$$
 $$=(a_1+k).{\overline \b},\;\; \overline \b\in  SpT_{2n}(\varpi_{t-1}), \s(\overline\b)=\overline\b.$$

 If $\overline\b=\emptyset\Leftrightarrow t=1$
 \begin{align}\label{b=emptythm+}
{\redu}_{t}^{-1}(\mathbf{b})=(1,2,\dots, l_1)(a_1+k)(a_1+k+1,\dots, a_1+k+1+(l-1-l_1))\in SST_{2n}(\varpi_l),
\end{align}
where
\begin{align}
&k+s\le l_1=
min\{a_1+k-2, l-1\}, 
\end{align}
and $$\a=(a_1, a_1+1,\dots ,a_1+{k-1})\sqcup (a_1+k)\sqcup (a_1+k+1,\dots, a_1+k+1+s-1)
 \in SST_{2n}(\varpi_{k+1+s}).$$

Then
\begin{align}
&{\redu}_{t}^{-1}(\mathbf{b})\bigcdot S^1={\redu}_{t}^{-1}(\mathbf{b})\bigcdot \a=\nonumber\\
 &=(1,2,\dots, l_1,a_1+k,a_1+k+1,\dots, a_1+k+1+l-1-l_1)\bigcdot \nonumber\\
 &\bigcdot(a_1, a_1+1,\dots ,a_1+{k-1})\sqcup (a_1+k)\sqcup (a_1+k+1,\dots, a_1+k+1+s-1)
 &\in SST_{2n}(\varpi_l+\varpi_{k+s+1}),\nonumber\\
 &  l\ge k+1+s,\; l_1\ge k+s.\nonumber
\end{align}

Otherwise, $\overline\b\neq\emptyset\Leftrightarrow 2\le t-1\in 2\Z$,

\begin{align} &\redu_{t}^{-1}(\mathbf{b})=(1,2,\dots, l_1)(a_1+k){\overline\b}(a_1+k-1+t+1,\dots, a_1+k-1+t+(l-t-l_1))\in
SST_{2n}(\varpi_l)\nonumber\\
&l_1\ge k+s,\; l\ge k+s+t\nonumber
\label{notemptythm+}
\end{align}
where $l_1$, $l_{t+1}$  are nonnegative even numbers defined by
\begin{align*}
k+s\le l_1=min\{a_1+k-1-t,l-t\},\quad
l_{t+1}=l-t-l_1,
\end{align*}


and \begin{align*}&{\redu}_{t}^{-1}(\mathbf{b})\bigcdot S^1={\redu}_{t}^{-1}(\mathbf{b})\bigcdot \a\nonumber\\
&(1,2,\dots, l_1)(a_1+k){\overline\b}(a_1+k-1+t+1,\dots, a_1+k-1+t+(l-t-l_1))\bigcdot\nonumber\\
&\bigcdot(a_1, a_1+1,\dots ,a_1+{k-1})\sqcup (a_1+k, a_1+k+1,\dots,a_1+k+t-1)\sqcup (a_1+k+t,\dots, a_1+k+t+s-1)\nonumber\\
&\in  SST_{2n}(\varpi_l+\varpi_{k+t+s}),\;
 l\ge k+t+s,\, l_1\ge k+s.
\end{align*}

$(c)$  if $ k,t, s\in [n]$, and $k\in 2\Z\Leftrightarrow a_1+k\notin 2\Z$, $s,t \notin 2\Z$ $\Rightarrow a_1+k+t-2\in 2\Z, \, a_1+k+t-1\notin 2\Z$,  it holds
$$\b=(a_1+k,a_1+k+1,\dots,a_1+k+t-2)(a_1+k+t-1)\in SpT_{2n}(\varpi_{t}), $$
 $$=\overline \b\bigcdot (a_1+k+t-1),\;\; \overline \b\in  SpT_{2n}(\varpi_{t-1}), \s(\overline\b)=\overline\b.$$
  Then
  \begin{align} &\redu_{t}^{-1}(\mathbf{b})=(1,2,\dots, l_1){\overline\b}(a_1+k+t-1)(a_1+k+t-1+2,\dots, a_1+k-1+t+1+(l-t-l_1))\in
SST_{2n}(\varpi_l)\nonumber\\
&l_1\ge k+s,\; l\ge k+s+t\nonumber
\end{align}
where $l_1$ is a nonnegative even number defined by
\begin{align*}
k+s\le l_1=min\{a_1+k-t,l-t\},\quad
\end{align*}
and
\begin{align*}&{\redu}_{t}^{-1}(\mathbf{b})\bigcdot \a\nonumber\\
&(1,2,\dots, l_1){\overline\b}(a_1+k+t-1)(a_1+k-1+t+2,\dots, a_1+k-1+t+1+(l-t-l_1))\bigcdot\nonumber\\
&\bigcdot(a_1, a_1+1,\dots ,a_1+{k-1})\sqcup (a_1+k, a_1+k+1,\dots,a_1+k+t-1)\sqcup (a_1+k+t,\dots, a_1+k+t+s-1)\nonumber\\
&\in  SST_{2n}(\varpi_l+\varpi_{k+t+s}),\;
 l\ge k+t+s,\, l_1\ge k+s.
\end{align*}

\medskip
 $(d)$  if $ k,t, s\in [n]$, and $t\in 2\Z$, $s,k \notin 2\Z$ $\Rightarrow a_1+k\in 2\Z$, $  \, a_1+k+t-1\notin 2\Z$,  it holds

 $$\b=(a_1+k)(a_1+k+1,\dots,a_1+k+t-2)(a_1+k+t-1)\in SpT_{2n}(\varpi_{t}), $$
 $$ =(a_1+k)\bigcdot\overline \b\bigcdot (a_1+k+t-1),\;\; \overline \b\in  SpT_{2n}(\varpi_{t-2}),\; \s(\overline \b)=\overline\b.$$

 Then
  \begin{align} &\redu_{t}^{-1}(\mathbf{b})=(1,2,\dots, l_1)(a_1+k){\overline\b}(a_1+k+t-1)(a_1+k+t-1+2,\dots, a_1+k-1+t+1+(l-t-l_1))\nonumber\\
  &\in
SST_{2n}(\varpi_l),
\;l_1\ge k+s,\; l\ge k+s+t\nonumber
\end{align}
where $l_1$ is a nonnegative even number defined by
\begin{align*}
k+s\le l_1=min\{a_1+k-t-1,l-t\},\quad
\end{align*}

and
\begin{align*}&{\redu}_{t}^{-1}(\mathbf{b})\bigcdot \a\nonumber\\
&(1,2,\dots, l_1)(a_1+k){\overline\b}(a_1+k+t-1)(a_1+k-1+t+2,\dots, a_1+k-1+t+1+(l-t-l_1))\bigcdot\nonumber\\
&\bigcdot(a_1, a_1+1,\dots ,a_1+{k-1})\sqcup (a_1+k, a_1+k+1,\dots,a_1+k+t-1)\sqcup (a_1+k+t,\dots, a_1+k+t+s-1)\nonumber\\
&\in  SST_{2n}(\varpi_l+\varpi_{k+t+s}),\;
 l\ge k+t+s,\, l_1\ge k+s.
\end{align*}

\end{enumerate}

Moreover,
\begin{align}
 &{\redu}_{t}^{-1}(\mathbf{b})\bigcdot (\underset{ (r_{k+1},\dots,r_{k+t})}\longleftarrow S)\in
 SST_{2n}(\varpi_l+\mu-\delta_{(r_{k+1},\dots,r_{k+t})}).
 \end{align}
 \qquad\qquad\qquad\qquad\qquad\qquad\qquad\qquad\qquad\qquad\qquad\qquad\qquad\qquad\qquad\qquad\qquad\qquad\qquad
 \qquad\qquad\qquad\qquad$\Box$
\end{thm}

\begin{ex}
\begin{enumerate}
\item Let $n=6$,
\begin{enumerate}
\item $k=2$, $t=4$, $s=0$, $k+t=6\in[6]\cap 2\Z$, $k+t\le l=8\le 2n$ where $t\le \min\{l=8,2n-l=4\}$ and $l-t=4\in 2\Z$. Let
$a_1=7$ and $a_1\ge k+t+s+1=7$,

\begin{align}\left(\b,\underset{ (2,3,5,6)}\longleftarrow S=\YT{0.15in}{}{
{7,8},
 {8,9},
 {9,10},
 {10},
 {11},
 {12},
}\right),\quad \b=(9,10,11,12)\in SpT_{2n}(\varpi_t), \; t=4
\end{align}
and $\a=(7,8,9,10,11,12)\in SpT_{2n}(\varpi_{k+t})$ such that $a_1=7\ge k+t+1\Leftrightarrow a_1+2\ge 2k+t+1$.

Then $\b=(a_1+2,a_1+3,a_1+4, a_1+5)\in SpT_{2n}(\varpi_t)$, and
$$\left(\underset{ (2,3,5,6)}\longleftarrow S\right)=\YT{0.15in}{}{
{7,8},
 {8},
 {9},
 {10},
}\in SST_{2n}(2,1,1,1)
$$

and \begin{align*}
&k\le l_1=9-5=4=min\{a_1+k-(t+1),l-t=4\},\quad l_1=l-t,\nonumber\\
&l_{t+1}=l-t-l_1=0\Leftrightarrow l_5=8-4-4=0.
\end{align*}
Hence

\begin{align}
{\redu}_{t}^{-1}(\mathbf{b})=(1,2,3, 4){\b},
\end{align}
and

\begin{align}{\redu}_{t}^{-1}(\mathbf{b})\bigcdot S^1=\YT{0.15in}{}{
{1,7,8},
 {2,8},
 {3,9},
 {4,10},
 {9},
 {10},
 {11},
 {12},
}\in SST_{2n}(\varpi_8+\varpi_4+\varpi_1)
\end{align}
\item $k=t=s=2\in [n]\cap 2\Z$, $k+t+s\le l=8\le 2n$ where $t\le \min\{l=8,2n-l=4\}$ and $l-t=6\in 2\Z$.
Let
$a_1=7$ and $a_1\ge k+t+s+1=7$,

\begin{align}\left(\b,\underset{ (2,3)}\longleftarrow S=\YT{0.15in}{}{
{7,8},
 {8,9},
 {9,10},
 {10},
 {11},
 {12},
}\right),\quad \b=(9,10)\in SpT_{2n}(\varpi_t), \; t=2
\end{align}

and $$\left(\underset{ (2,3)}\longleftarrow S\right)=\YT{0.15in}{}{
{7,8},
 {8},
 {9},
 {10},
 {11},
 {12},
}\in SST_{2n}(\varpi_6+\varpi_1)
$$
and \begin{align*}
&k\le l_1=6=min\{a_1+k-(t+1),l-t\}=\{7+2-2-1,8-2\},\; l_1=l-t\nonumber\\
&l_{t+1}=l-t-l_1=0\Leftrightarrow l_5=8-2-6=0.
\end{align*}
Hence

\begin{align}
{\redu}_{t}^{-1}(\mathbf{b})=(1,2,3, 4,5,6){\b},
\end{align}
and

\begin{align}{\redu}_{t}^{-1}(\mathbf{b})\bigcdot \left(\underset{ (2,3)}\longleftarrow S\right)=\YT{0.15in}{}{
{1,7,8},
 {2,8},
 {3,9},
 {4,10},
 {5,11},
 {6,12},
 {9},
 {10},
}\in SST_{2n}(\varpi_8+\varpi_6+\varpi_1)
\end{align}
\item  $k=t=2=s\in [n]\cap 2\Z$, $k+t+s\le l=10\le 2n$ where $2=t\le \min\{l=10,2n-l=2\}$ and $l-t=8\in 2\Z$.
Let
$a_1=7$ and $a_1\ge k+t+s+1=7$,

\begin{align}\left(\b,\underset{ (2,3)}\longleftarrow S=\YT{0.15in}{}{
{7,8},
 {8,9},
 {9,10},
 {10},
 {11},
 {12},
}\right),\quad \b=(9,10)\in SpT_{2n}(\varpi_t), \; t=2
\end{align}

and $$\left(\underset{ (2,3)}\longleftarrow S\right)=\YT{0.15in}{}{
{7,8},
 {8},
 {9},
 {10},
 {11},
 {12},
}\in SST_{2n}(\varpi_6+\varpi_1)
$$
and \begin{align*}
&k\le l_1=6=min\{a_1+k-(t+1),l-t\}=\{7+2-2-1,10-2\},l_1=6<l-t\nonumber\\
&l_{t+1}=l-t-l_1\Leftrightarrow l_5=10-2-6=2.
\end{align*}
Hence

\begin{align}
{\redu}_{t}^{-1}(\mathbf{b})=(1,2,3, 4,5,6){\b}(11,12)\in SST_{2n}(\varpi_{10}),
\end{align}
and

\begin{align}{\redu}_{t}^{-1}(\mathbf{b})\bigcdot \left(\underset{ (2,3)}\longleftarrow S\right)=\YT{0.15in}{}{
{1,7,8},
 {2,8},
 {3,9},
 {4,10},
 {5,11},
 {6,12},
 {9},
 {10},
 {11},
 {12},
}\in SST_{2n}(\varpi_{10}+\varpi_6+\varpi_1)
\end{align}
\end{enumerate}
\item Let $n=7$,
\begin{enumerate}
\item
$k=t=2=s\in [n]\cap 2\Z$, $k+t+s\le l=12\le 2n$ where $2=t\le \min\{l=12,2n-l=2\}$ and $l-t=12-2\in 2\Z$.
Let
$a_1=7$ and $a_1\ge k+t+s+1=7$,

\begin{align}\left(\b,\underset{ (2,3)}\longleftarrow S=\YT{0.15in}{}{
{7,8},
 {8,9},
 {9,10},
 {10},
 {11},
 {12},
}\right),\quad \b=(9,10)\in SpT_{2n}(\varpi_t), \; t=2
\end{align}

and $$\left(\underset{ (2,3)}\longleftarrow S\right)=\YT{0.15in}{}{
{7,8},
 {8},
 {9},
 {10},
 {11},
 {12},
}\in SST_{2n}(\varpi_6+\varpi_1)
$$
and \begin{align*}
&k\le l_1=6=min\{a_1+k-(t+1),l-t\}=\{7+2-2-1,12-2\},\; l_1<l-t\nonumber\\
&l_{t+1}=l-t-l_1\Leftrightarrow l_5=12-2-6=4.
\end{align*}
Hence

\begin{align}
{\redu}_{t}^{-1}(\mathbf{b})=(1,2,3, 4,5,6){\b}(11,12,13,14)\in SST_{2n}(\varpi_{12}),
\end{align}
and

\begin{align}{\redu}_{t}^{-1}(\mathbf{b})\bigcdot \left(\underset{ (2,3)}\longleftarrow S\right)=\YT{0.15in}{}{
{1,7,8},
 {2,8},
 {3,9},
 {4,10},
 {5,11},
 {6,12},
 {9},
 {10},
 {11},
 {12},
 {13},
 {14},
}\in SST_{2n}(\varpi_{12}+\varpi_6+\varpi_1)
\end{align}
\end{enumerate}
\item Let $n=9$,
\begin{enumerate}

\item
$k=t=2=s\in [n]\cap 2\Z$, $k+t+s\le l=16\le 2n$ where $2=t\le \min\{l=16,2n-l=2\}$ and $l-t=16-2\in 2\Z$.
Let
$a_1=7$ and $a_1\ge k+t+s+1=7$,
and \begin{align*}
&k\le l_1=6=min\{a_1+k-(t+1),l-t\}=\{7+2-2-1,16-2\},\; l_1<l-t\nonumber\\
&l_{t+1}=l-t-l_1\Leftrightarrow l_5=16-2-6=8.
\end{align*}
Hence

\begin{align}
{\redu}_{t}^{-1}(\mathbf{b})=(1,2,3, 4,5,6){\b}(11,12,13,14,15,16,17,18)\in SST_{2n}(\varpi_{12}),
\end{align}
where the interval ${\redu}_{t}^{-1}(\mathbf{b})=(1,2,3, 4,5,6){\b}(11,12,13,14,15,16,17,18)$ is not symplectic,
and

\begin{align}{\redu}_{t}^{-1}(\mathbf{b})\bigcdot \left(\underset{ (2,3)}\longleftarrow S\right)=\YT{0.15in}{}{
{1,7,8},
 {2,8},
 {3,9},
 {4,10},
 {5,11},
 {6,12},
 {9},
 {10},
 {11},
 {12},
 {13},
 {14},
 {15},
 {16},
 {17},
 {18}
}\in SST_{2n}(\varpi_{16}+\varpi_6+\varpi_1)
\end{align}

\end{enumerate}
\end{enumerate}
\end{ex}
\subsubsection{The bumped column is of the form $x\B y$ such that $\s(\B)=\B$}\label{bumpedxintervaly}

\medskip
The following theorem is a consequence of Theorem \ref{thm:xAy} with $p_0=p_1=1$.

\begin{thm}\label{thm:bumpedxintervaly}
Let $k+t+s\in[n]\cap 2\Z$, for some $ k,s, t\in 2\Z$ and  
$a_1\notin 2\Z$. Let $1\le u\le k$ 
  and $1\le v\le s$. Let $S\in SpT_{2n}(\mu)$ be a symplectic tableau with first  column  the symplectic interval $\a$ in $ SpT_{2n}(\varpi_{k+t+s})$ such that $\s(\a)=\a$,
\begin{align}
&\mathbf{a}
= (a_1, a_1+1,\dots ,\mathsf{a_1+u-1},\dots,a_1+{k-1})\sqcup \nonumber\\
&\sqcup (a_1+k, a_1+k+1,\dots,a_1+k+t-1)\sqcup \nonumber\\
&\sqcup (a_1+k+t,a_1+k+t+1\dots,\mathsf{a_1+k+t+v-1}, \dots,a_1+k+t+s-1)
\end{align}
where (necessarily)
$a_1\ge k+t+s+1$.
Let $r_u<r_{k+1}<\dots<r_{k+t}<r_{k+t+v}$ be $t+2$ row coordinates of   $S$  such that
$\mu-\delta_{(r_u,r_{k+1},\dots,r_{k+t},r_{k+t+v})}$ is a partition, and
$$(\b, \underset{ (r_u,r_{k+1},\dots,r_{k+t},r_{k+t+v})}\longleftarrow S),$$
with
$$\b=(\mathsf{a_1+u-1})\sqcup(a_1+k, a_1+k+1,\dots,a_1+k+t-1)\sqcup(\mathsf{a_1+k+t-1+v})\in SpT_{2n}(\varpi_{t+2}),
$$   and  $$(\underset{(r_u,r_{k+1},\dots,r_{k+t},r_{k+t+v})}\longleftarrow S)\in
SpT_{2n}(\mu-\delta_{(r_u,r_{k+1},\dots,r_{k+t},r_{k+t+v})})$$
has first column $S^1\in SpT_{2n}(\varpi_{k+t+s}).$
Then
\begin{enumerate}
\item
$S^1=\a\in SpT_{2n}(\varpi_{k+t+s})$ 
and 
  $ a_1+k\notin 2\Z$, $ a_1+k+t-1\in 2\Z$.
  \medskip
  \item $\b=(\mathsf{a_1+u-1})\bigcdot \overline\b  \bigcdot (\mathsf{a_1+k+t+v-1})$, with $\overline\b=(a_1+k, a_1+k+1,\dots,a_1+k+t-1)$,  
  such that $\s(\overline\b)=\overline\b$ and $\s((\mathsf{a_1+u-1}).(\mathsf{a_1+k+t+v-1}))\cap \b=\emptyset$.

 \medskip
   
   \item 
$a_1+k\ge 2k+t+s+1\Rightarrow a_1+k-t-2-1\ge 2k+t+s+1-t-2-1\ge 2k+s-2\ge k+s.$

\medskip
\item For $k+t+s\le l\le 2n$, and $k,t,s\in[n]\cap 2\Z$
such that
 $ t+2\le min\{l,2n-l\}$ and $l-t-2\in 2\Z$, considering the even numbers
$$\tilde l_1=min\{a_1+k-t-2-1, l-t-2\}\ge k+s-2,$$
where
$$\tilde l_1\ge\begin{cases} k,& s=2\in 2\Z,\\
k+2,& s\ge 4, \;s\in 2\Z,
\end{cases}
$$
and
\begin{align}l_1=\begin{cases}min\{a_1+u-3,\tilde l_1\}\ge k+s-2\ge k,& \mbox{ if  } a_1+u-1\in 2\Z\Leftrightarrow u\in 2\Z\\ 
min\{a_1+u-2, \tilde l_1\}\ge k+s-2\ge k, &  \mbox{ if  } a_1+u-1\notin 2\Z\Leftrightarrow u\notin 2\Z
\end{cases}, 
\end{align}
where 
$$ l_1\ge\begin{cases} k,& s=2\in 2\Z,\\
k+2,& s\ge 4, \;s\in 2\Z,
\end{cases}$$
one defines
\begin{align}T_1=\begin{cases}(\mathsf{a_1+{u-1}}+1,\dots, \mathsf{a_1+{u-1}}+l_2),& \mbox{ if } \mathsf{a_1+{u-1}}\in 2\Z\\
(\mathsf{a_1+{u-1}}+2,\dots,\mathsf{a_1+{u-1}}+1+l_2),& \mbox{ if } \mathsf{a_1+{u-1}}\notin 2\Z,\quad \; l_2=\tilde l_1-l_1\ge 0
\end{cases},
\end{align}
and 
$$(1,2,\dots, l_1)\sqcup (\mathsf{a_1+{u-1}})\sqcup T_1\sqcup \overline\b,$$
where
$$\ell(1,2,\dots, l_1)\sqcup (\mathsf{a_1+{u-1}}))= l_1+1\ge k+1\ge u,$$
$$\ell(1,2,\dots, l_1)\sqcup (\mathsf{a_1+{u-1}})\sqcup T_1\sqcup \overline\b)=l_1+l_2+t+1\ge 
\begin{cases}k+t+l_2+1,& \mbox{ if $s=2$}\\
k+t+l_2+3,& \mbox{ if $s\ge 4$.}
\end{cases}
$$

Moreover, considering the even numbers
$$\tilde l_2=l-t-2-l_1-l_2\ge k+t+s-2-t-l_1-l_2=k+s-2-l_1-l_2,$$
\begin{align}l_{t+2}=\begin{cases}min\{v-2,\tilde l_2\},&v\in 2\Z\\
min\{v-1,\tilde l_2\}, & v\notin 2\Z,\\
\end{cases}\Rightarrow l_{t+2}=0, \mbox{ if $v=1,2$ or $s=2$}
\end{align}
and $$l_{t+3}=\tilde l_2-l_{t+2}=l-t-2-l_1-l_2-l_{t+2},$$

one defines
$$Y_0=(a_1+k+(t-1)+1,\dots,a_1+k+(t-1)+l_{t+2}),$$

and, for $y_1=\mathsf{a_1+k+t+v-1}\Rightarrow y_1+2=\mathsf{a_1+k+t+v-1}+2=a_1+k+t+v+1$,
\begin{align}Y_1=\begin{cases}(y_1+1,\dots, y_1+l_{t+1+2})& \mbox{ if $y_1\in 2\Z\Leftrightarrow v\in 2\Z$ }\\
(y_1+1+1,\dots, y_1+1+l_{t+1+2})& \mbox{ if $y_1\notin 2\Z\Leftrightarrow v\notin 2\Z$ },
\end{cases}
\end{align}
where
$$\ell(1,2,\dots, l_1)\sqcup (\mathsf{a_1+{u-1}})\sqcup T_1\sqcup \overline\b\sqcup Y_0\sqcup(y_1))=l_1+l_2+t+l_{t+2}+1+1$$
$$\ge
\begin{cases}k+t+l_2+1+1,& \mbox{ if $s=2$}\\
k+t+l_2+l_{t+2}+4,& \mbox{ if $s\ge 4$.}
\end{cases}
$$

\bigskip
Therefore, it holds 
\begin{align} \label{redubthmxby}
&{\redu}_{t}^{-1}(\mathbf{b})=(1,2,\dots, l_1)\sqcup(\mathsf{a_1+{u-1}})\sqcup T_1\sqcup {\overline\b}\sqcup Y_0\sqcup (y_1)Y_1
\in SST_{2n}(\varpi_l),\nonumber\\
& \;\; l=l_1+1+l_2+t+l_{t+2}+1+l_{t+3},\nonumber
\end{align}
such that 
\begin{align}
&{\redu}_{t}^{-1}(\mathbf{b})(k+1,1)  \le {\redu}_{t}^{-1}(\mathbf{b})(l_1+1,1)=\mathsf{a_1+{u-1}} \le a_1+k-1<S^1(k+1,1)=a_1+k, \; l_1\ge k,\nonumber
\end{align}
\begin{align}
&{\redu}_{t}^{-1}(\mathbf{b})(k+t+1,1)\le{\redu}_{t}^{-1}(\mathbf{b})(l_1+l_2+t+1,1)=a_1+k+t-1<S^1(k+t+1,1)=a_1+k+t \nonumber 
\end{align}

\underline{For $s=2$}, \underline{then  either $v=1$ or $v=2$}, and $l_{t+2}=0$, $Y_0=\emptyset$ and $\ell(S^1)=k+t+2$. 

If $v=1$,  $y_1=\mathsf{a_1+k+t+v-1}= a_1+k+t$,  and 
\begin{align}\oo
&{\redu}_{t}^{-1}({\b})(k+t+2,1)\le\nonumber\\
&\le {\redu}_{t}^{-1}(\mathbf{b})(l_1+l_2+t+2,1)=a_1+k+t=y_1<S^1(k+t+2,1)=a_1+k+t+1. \nonumber
\end{align}
If $v=2$, $y_1=\mathsf{a_1+k+t+v-1}= a_1+k+t+1$, and 
\begin{align}\oo
&{\redu}_{t}^{-1}({\b})(k+t+2,1)\le{\redu}_{t}^{-1}(\mathbf{b})(l_1+l_2+t+2,1)=y_1=S^1(k+t+2,1)=a_1+k+t+1. \nonumber
\end{align}

\underline{ For $s\ge 4$},  and $v=1$ or $v=2$, one has $l_{t+2}=0$, $Y_0=\emptyset$ and $l\ge \ell(S^1)=k+t+s\ge k+t+4$.
 
 If $v=1$, $y_1=\mathsf{a_1+k+t+v-1}= a_1+k+t$
 \begin{align}\oo
&{\redu}_{t}^{-1}({\b})(k+t+3,1)\le\nonumber\\
&\le{\redu}_{t}^{-1}(\mathbf{b})(l_1+l_2+t+3,1)=a_1+k+t+2=y_1+2=S^1(k+t+3,1)=a_1+k+t+2. \nonumber
\end{align}

If $v=2$, $y_1=\mathsf{a_1+k+t+v-1}= a_1+k+t+1$,
 \begin{align}\oo
&{\redu}_{t}^{-1}({\b})(k+t+3,1)\le\nonumber\\
&\le{\redu}_{t}^{-1}(\mathbf{b})(l_1+l_2+t+3,1)=a_1+k+t+2=y_1+1=S^1(k+t+3,1)=a_1+k+t+2. \nonumber
\end{align}
 
 \underline{ For $s\ge 4$},  and \underline{$v\ge 3$}, one has  $l_1\ge k+2$,  and 
  \begin{align}& 2\Z\not\ni\ell( (1,2,\dots, l_1)\sqcup(\mathsf{a_1+{u-1}})\sqcup T_1\sqcup {\overline\b})=\tilde l_1+t+1=l_1+l_2+t+1\ge k+s-1+l_2+t\nonumber\\
&\ge  k+t+v-1 \mbox{ if $v\in 2\Z$ and $k+t+v$ if  $v\notin 2\Z$},\nonumber
\end{align}
and \begin{align}& 2\Z\not\ni\ell( (1,2,\dots, l_1)\sqcup(\mathsf{a_1+{u-1}})\sqcup T_1\sqcup {\overline\b}\sqcup Y_0)=\tilde l_1+t+1+l_{t+2}=l_1+l_2+t+1+l_{t+2}\ge \nonumber\\
&\ge k+s-1+l_2+t++l_{t+2}\nonumber\\
&\ge  k+t+v-1 +l_{t+2}\mbox{ if $v\in 2\Z$ and $k+t+v$ if  $v\notin 2\Z$}.\nonumber
\end{align}

Then
  \begin{align}
&{\redu}_{t}^{-1}(\mathbf{b})(k+t+2,1)< {\redu}_{t}^{-1}(\mathbf{b})(k+t+4,1)\le \nonumber\\
&\le {\redu}_{t}^{-1}(\mathbf{b})(l_1+l_2+t+2,1)=a_1+k+t<S^1(k+t+2,1)=a_1+k+t+1 \nonumber 
\end{align}
\begin{align}
&\vdots\nonumber\\
&{\redu}_{t}^{-1}(\mathbf{b})(k+t+l_{t+2}+1,1)<\nonumber\\
&<{\redu}_{t}^{-1}(\mathbf{b})(l_1+l_2+t+l_{t+2}+1,1)=a_1+k+(t-1)+l_{t+2}<\nonumber\\
&< S^1(k+t+l_{t+2}+1,1)=a_1+k+t+l_{t+2}.\nonumber 
\end{align}
Moreover $(1,2,\dots, l_1)\sqcup(\mathsf{a_1+{u-1}})\sqcup T_1\sqcup {\overline\b}\sqcup Y_0\sqcup (y_1)$ is a column,
$$a_1+k+(t-1)+l_{t+2}<y_1=\mathsf{a_1+k+t+v-1},$$ 
\begin{align}&{\redu}_{t}^{-1}(\mathbf{b})(k+t+v,1)<\nonumber\\
&<{\redu}_{t}^{-1}(\mathbf{b})(l_1+l_2+t+1+l_{t+2}+1,1)=y_1
=\mathsf{a_1+k+t+v-1}=S^1(k+t+v,1) \nonumber
\end{align}
and
\begin{align}&{\redu}_{t}^{-1}(\mathbf{b})(l_1+1+l_2+t+l_{t+2}+1+1,1)=\nonumber\\
&=\begin{cases}
y_1+1=a_1+k+t+v-1+1=S^1(k+t+v+1,1)& \mbox{ if  $v\in 2\Z$ }\\
y_1+2=a_1+k+t+v-1+2=S^1(k+t+v+2,1)& \mbox{ if  $v\notin 2\Z$ }
\end{cases}
\end{align}
where
\begin{align}
&{\redu}_{t}^{-1}(\mathbf{b})(k+t+v+1,1)<{\redu}_{t}^{-1}(\mathbf{b})(l_1+1+l_2+t+l_{t+2}+1+1,1) \mbox{ if  $v\in 2\Z$ }\nonumber\\
&{\redu}_{t}^{-1}(\mathbf{b})(k+t+v+2,1)<{\redu}_{t}^{-1}(\mathbf{b})(l_1+1+l_2+t+l_{t+2}+1+1,1) \mbox{ if  $v\notin 2\Z$ } \nonumber
\end{align}

That is, the  concatenation of $\redu_{t}^{-1}(\b)$ and $S^1$ 
 is a semi-standard tableau in $ SST_{2n}(\varpi_l+\varpi_{k+t+s})$, with $l\ge k+t+s$,
\begin{align}&{\redu}_{t}^{-1}(\mathbf{b})\bigcdot S^1={\redu}_{t}^{-1}(\mathbf{b})\bigcdot\a \nonumber\\
 &=(1,2,\dots, l_1)\sqcup(\mathsf{a_1+{u-1}})\sqcup T_1\sqcup {\overline\b}\sqcup Y_0\sqcup (y_1)Y_1\bigcdot\nonumber\\
 &\bigcdot(a_1, a_1+1,\dots ,a_1+{k-1})\sqcup (a_1+k, a_1+k+1,\dots,a_1+k+t-1)\sqcup (a_1+k+t,\dots, a_1+k+t+s-1)\nonumber\\
 &\in SST_{2n}(\varpi_l+\varpi_{k+t+s}),\nonumber
 \end{align}
and
 \begin{align}
 &{\redu}_{t}^{-1}(\mathbf{b})\bigcdot (\underset{(r_u,r_{k+1},\dots,r_{k+t},r_{k+t+v})}\longleftarrow S)\in
SpT_{2n}(\varpi_l+\mu-\delta_{(r_u,r_{k+1},\dots,r_{k+t},r_{k+t+v})})
 \end{align}
 \end{enumerate}
 $\qquad\qquad\qquad\qquad\qquad\qquad\qquad\qquad\qquad\qquad\qquad\qquad\qquad\qquad\qquad\qquad\qquad\qquad\qquad\qquad\qquad\qquad\qquad\Box$
\end{thm}

\subsection{Bumped symplectic intervals from semistandard tableaux}\label{bumpedinterval}
{
We consider a symplectic bumped interval $\b$, not necessarily fixed by the parity involution, without imposing that the tableau from where it is bumped out is symplectic.}

\begin{thm}\label{thm:reducbumpintervalgeneral} 

Let $t\in[n]$,
 and  $S\in SST_{2n}(\mu)$ with first column $\mathbf{a} = (a_1, \dots ,a_{k})\sqcup \b \sqcup (a_{k+t+1}, \dots,a_{k+t+s})\in
 SST_{2n}(\varpi_{k+t+s})$ such that $\b=(b_1,b_1+1,\dots,b_1+{t-1})$ , $b_1\notin 2\Z$ and $ b_1\ge max\{a_k+1,k+t+1\}$.
Let $r_{k+1}<\dots<r_{k+t}$ be $t$ row coordinates of   $S$ such that
$\mu-\delta_{(r_{k+1},\dots,r_{k+t})}$ is a partition, and
$$(\b, \underset{ (r_{k+1},\dots,r_{k+t})}\longleftarrow S),$$
where 
$$\b=(b_1, b_1+1,\dots,b_1+t-1)\in SST_{2n}(\varpi_t)$$   and  $$(\underset{(r_{k+1},\dots,r_{k+t})}\longleftarrow S)\in
SST_{2n}(\mu-\delta_{(r_{k+1},\dots,r_{k+t})})$$
has first column $S^1=(a_1, a_2,\dots ,a_k)\sqcup (x_k,x_{k+1},\dots, x_{k+t-1})\sqcup (a_{k+t+1}, \dots,a_{k+t+s})\in SST_{2n}(\varpi_{k+t+s}).$
Then
\begin{enumerate}
\item $\b=(b_1, b_1+1,\dots,b_1+t-1)\in SpT_{2n}(\varpi_t)$ and $S^1\in SST_{2n}(\varpi_{k+t+s})$ are such that

 $$(x_k,x_{k+1}\dots, x_{k+t-1})\ge \b=(b_1, b_1+1,\dots, b_1+t-1).$$
\item For $k+t+s\le l\le 2n$, 
such that
 $ t\le min\{l,2n-l\}$ and $l-t\in 2\Z$,

\medskip
$(a)$ if $t\in 2\Z$ and  $b_1\notin 2\Z$ $\Rightarrow b_1+t-1\in 2\Z$, and it holds $\s(\b)=\b$ and

\begin{align} 
{\redu}_{t}^{-1}(\mathbf{b})=(1,2,\dots, l_1){\b}(b_1+t-1+1,\dots, b_1+t-1+(l-t-l_1))\in SST_{2n}(\varpi_l),
\end{align}
where
and $l_1$, $l_{t+1}$  are nonnegative even numbers defined by
\begin{align*}
k\le l_1=min\{b_1-t-1,l-t\},\quad
l_{t+1}=l-t-l_1,
\end{align*}

and
 \begin{align}&{\redu}_{t}^{-1}(\mathbf{b})\bigcdot S^1=\nonumber\\
 &=(1,2,\dots, l_1)\sqcup(b_1, \dots,b_1+t-1)\sqcup(b_1+t-1+1,\dots, b_1+t-1+(l-t-l_1))\bigcdot\nonumber\\
 &\bigcdot(a_1, \dots ,a_{1+{k-1}})\sqcup(x_k,x_{k+1}\dots, x_{k+t-1})\sqcup(a_{k+t+1}, \dots,a_{k+t+s})\in SST_{2n}(\varpi_l+\varpi_{k+t+s}),\nonumber\\
 & \mbox{ with $ l_1\ge k,$ and $l\ge  k+t+s$}.
 \end{align}
\medskip

$(b)$ if $t\notin 2\Z$ and $b_1\in 2\Z$  $\Rightarrow b_1+1\notin 2\Z,          \, b_1+t-1\in 2\Z$, it holds

 $$\b=(b_1)\sqcup (b_1+1,\dots,b_1+t-1)=(b_1)\overline \b\in SpT_{2n}(\varpi_{t}),\; b_1\in 2\Z$$
 $$ \overline \b\in  SpT_{2n}(\varpi_{t-1}),\;\s(\overline \b)=\overline \b.$$

Then
\begin{align} &\redu_{t}^{-1}(\mathbf{b})=(1,2,\dots, l_1)(b _1){\overline\b}(b_1+t-1+1,\dots, b_1+t-1+(l-(t-1)-1-l_1))\in  SST_{2n}(\varpi_l,\\
&l_1\ge k,\; l\ge k+t+s
\end{align}

where $k\le l_1=min\{b_1+1-(t-1+3), l-(t-1+1)\}$.
and \begin{align*}&{\redu}_{t}^{-1}(\mathbf{b})\bigcdot S^1=\nonumber\\
&{\redu}_{t}^{-1}(\mathbf{b})\bigcdot (a_1, a_1+1,\dots ,a_1+{k-1})\sqcup (x_k,x_{k+1}\dots, x_{k+t-1})\\
&\sqcup
(a_{k+t+1}, \dots,a_{k+t+s})\in  SST_{2n}(\varpi_l+\varpi_{k+t+s}),\; l\ge k+t+s,\, l_1\ge k.
\end{align*}

$(c)$ if $t\notin 2\Z$ and $b_1\notin 2\Z$  $\Rightarrow          \, b_1+t-2\in 2\Z$, it holds

 $$\b=(b_1,\dots,b_1+t-2)(b_1+t-1)=\overline \b(b_1+t-1)\in SpT_{2n}(\varpi_{t}),\; b_1+t-1\notin 2\Z$$
 $$ \overline \b\in  SpT_{2n}(\varpi_{t-1}),\;\s(\overline \b)=\overline \b.$$
 Then

 \begin{align} &\redu_{t}^{-1}(\mathbf{b})=(1,2,\dots, l_1){\overline\b}(b_1+t-1)(b_1+t+1,\dots, b_1+t+1+(l-(t-1)-1-l_1))\in  SST_{2n}(\varpi_l,\\
&l_1\ge k+1,\; l\ge k+t+s
\end{align}
  where
 $k+1\le l_1=min\{b_1-(t-1+1),l-t)\}$

 and \begin{align*}&{\redu}_{t}^{-1}(\mathbf{b})\bigcdot S^1=\nonumber\\
&{\redu}_{t}^{-1}(\mathbf{b})\bigcdot (a_1, a_1+1,\dots ,a_1+{k-1})\sqcup (x_k,x_{k+1}\dots, x_{k+t-1})\\
&\sqcup
(a_{k+t+1}, \dots,a_{k+t+s})\in  SST_{2n}(\varpi_l+\varpi_{k+t+s}),\; l\ge k+t+s,\, l_1\ge k+1.
\end{align*}

$(d)$ if $t\in 2\Z$ and $b_1\in 2\Z$  $\Rightarrow          \, b_1+t-2\in 2\Z$, it holds

 $$\b=(b_1)(b_1+1,\dots,b_1+t-2)(b_1+t-1)=(b_1)\overline \b(b_1+t-1)\in SpT_{2n}(\varpi_{t}),\;b_1\in 2\Z,\; b_1+t-1\notin 2\Z$$
 $$ \overline \b\in  SpT_{2n}(\varpi_{t-2}),\;\s(\overline \b)=\overline \b.$$

 Then
\begin{align} &\redu_{t}^{-1}(\mathbf{b})=(1,2,\dots, l_1)(b_1){\overline\b}(b_1+t-1)(b_1+t+1,\dots, b_1+t+1+(l-(t-2)-2-l_1))\nonumber\\
 &\in  SST_{2n}(\varpi_l,
\;l_1\ge k+1,\; l\ge k+t+s
\end{align}
  where
 $k+1\le l_1=min\{b_1+1-(t-2+3),l-t)\}$
 and \begin{align*}&{\redu}_{t}^{-1}(\mathbf{b})\bigcdot S^1=\nonumber\\
&{\redu}_{t}^{-1}(\mathbf{b})\bigcdot (a_1, a_1+1,\dots ,a_1+{k-1})\sqcup (x_k,x_{k+1}\dots, x_{k+t-1})\\
&\sqcup
(a_{k+t+1}, \dots,a_{k+t+s})\in  SST_{2n}(\varpi_l+\varpi_{k+t+s}),\; l\ge k+t+s,\, l_1\ge k+1.
\end{align*}
 \end{enumerate}
 \end{thm}
 \begin{proof} $(1)$ Follows from the condition $b_1\ge k+t+1$ and  properties of reverse Schenested insertion in Section \ref{subsec:insert}.

 $(2)$ Follows from Lemma \ref{lem:reductionsurj}, and corollaries \ref{cor:xA}, \ref{cor:Ay} and \ref{cor:xAy}.
 \end{proof}

\bibliography{sample17}
\bibliographystyle{alpha}
\end{document}